\documentclass[reqno,11pt]{amsart}

\usepackage{amssymb,amsmath}
\usepackage{amscd}

\theoremstyle{plain}
  \newtheorem{theorem}{Theorem}[section]
  \newtheorem{corollary}[theorem]{Corollary}
  \newtheorem{lemma}[theorem]{Lemma}
  \newtheorem{proposition}[theorem]{Proposition}
  \newtheorem{conjecture}[theorem]{Conjecture}
  \newtheorem{observation}[theorem]{Observation}
\theoremstyle{definition}
  \newtheorem{definition}[theorem]{Definition}

  \newtheorem{remark}[theorem]{Remark}

\DeclareMathOperator{\rank}{rank}

\DeclareMathOperator{\ind}{ind}

\DeclareMathOperator{\Int}{Int}
\DeclareMathOperator{\grad}{grad}

\DeclareMathOperator{\id}{id}

\DeclareMathOperator{\FT}{FT}
\newcommand\ii{\mathbf i}
\renewcommand\epsilon\varepsilon
\newcommand\rrangle{\rangle\!\rangle}
\newcommand\llangle{\langle\!\langle}
\renewcommand\Re{\operatorname{Re}}
\newcommand\LO{\operatorname{O}}

\newcommand\supp{\operatorname{supp}}
\newcommand\Spec{\operatorname{Spec}}

\newcommand\vvol{\operatorname{vol}}

\newcommand\tr{\operatorname{tr}}

\newcommand\eend{\operatorname{end}}

\newcommand\sdet{\operatorname{sdet}}

\newcommand{\N}{\mathbb N}

\newcommand{\R}{\mathbb R}
\newcommand{\C}{\mathbb C}

\newcommand{\sm}{\text{\rm sm}}

\newcommand{\la}{\text{\rm la}}

\def\co{\colon\thinspace}

\begin{document}

\title{Complex valued Ray--Singer torsion II}

\author{Dan Burghelea}

\address{Department of Mathematics,
         The Ohio State University,\\ 
         231 West 18th Avenue,
         Columbus, OH 43210, USA.}

\email{burghele@mps.ohio-state.edu}

\author{Stefan Haller}

\address{Department of Mathematics,
         University of Vienna,\\
         Nordbergstra{\ss}e 15,
         A-1090 Vienna, Austria.} 

\email{stefan.haller@univie.ac.at}

\thanks{Part of this work was done while the second author enjoyed the
        worm hospitality of the Ohio State University.}




%

\begin{abstract}
In this paper we extend Witten--Helffer--Sj\"ostrand theory from selfadjoint Laplacians based on fiber wise
Hermitian structures, to non-selfadjoint Laplacians based on fiber wise non-degenerate symmetric bilinear forms.
As an application we verify, up to sign, the conjecture about the comparison of the Milnor--Turaev torsion with the complex valued 
analytic torsion, for odd dimensional manifolds. This is done along  the lines of 
Burghelea, Friedlander and Kappeler's  proof of  the Cheeger--M\"uller theorem.
\end{abstract}

\maketitle
 
\setcounter{tocdepth}{1}
\tableofcontents

\section{Introduction}\label{S:intro}

Let $(M,g)$ be a closed connected smooth Riemannian manifold of dimension $n$, and suppose
$E$ is a flat complex vector bundle over $M$, equipped with a (not necessarily parallel) fiber wise 
non-degenerate symmetric bilinear form $b$. The flat connection of $E$ will be denoted by $\nabla^E$. 
Let $\Omega(M;E)$ denote the deRham complex of $E$-valued differential forms on $M$,
and write $d_E$ for the deRham differential.

The Riemannian metric $g$ and the bilinear form $b$ provide a fiber wise 
non-degenerate symmetric bilinear form on the complex vector bundle
$\Lambda^*T^*M\otimes E$ which will be denoted by $b_g$. If $\vvol_g$ denotes the volume density
associated with $g$, then 
$$
\beta(v,w):=\int_Mb_g(v,w)\vvol_g \qquad v,w\in\Omega(M;E),
$$
defines a symmetric bilinear form on $\Omega(M;E)$.
Let $d^\sharp_{E,g,b}$ denote the formal transposed of $d_E$ with respect to $\beta$,
$$
\beta(d_Ev,w)=\beta(v,d^\sharp_{E,g,b}w),\qquad
v,w\in\Omega(M;E).
$$
The Laplace--Beltrami operator 
\begin{equation}\label{E:Y2}
\Delta_{E,g,b}:=(d_E+d^\sharp_{E,g,b})^2=d_Ed^\sharp_{E,g,b}+d^\sharp_{E,g,b}d_E
\end{equation}
is not necessarily selfadjoint.

In \cite{BH06b} the following complex valued analogue of the square of the Ray--Singer torsion \cite{RS71} was introduced and studied:
\begin{equation}\label{E:crs}
\prod_q\bigl({\det}'(\Delta_{E,g,b,q})\bigr)^{(-1)^qq}\in\C^\times:=\C\setminus\{0\}.
\end{equation}
Here ${\det}'(\Delta_{E,g,b,q})\in\C^\times$ denotes the zeta regularized product of all non-vanishing eigen values of the 
Laplacian acting on $\Omega^q(M;E)$. In the definition of the determinant one can use any non-zero Agmon angle, the resulting ${\det}'(\Delta_{E,g,b,q})$ will
be independent of this choice.

Let $X$ be a Morse--Smale vector field\footnote{That is, with respect to some Riemannian metric, $X$ is the negative 
gradient of a Morse function, and satisfies the Smale transversality condition, i.e.\ stable and unstable manifolds 
intersect transversally.} 
on $M$, and write $\mathcal X$ for the set of critical points, i.e.\ zeros of $X$.
Denote by $C(X;E)$ the associated Morse complex. 
The fiber wise symmetric bilinear form $b$ induces a non-degenerate bilinear form $b_{\mathcal X}$ on $C(X;E)$.
Recall that the integration homomorphism
\begin{equation}\label{E:Int}
\Int_{E,X}:\Omega(M;E)\to C(X;E)
\end{equation}
induces an isomorphism in cohomology, see \cite{S93}.
Let $\Omega_{g,b}(M;E)(0)$ denote the (generalized) zero eigen space of $\Delta_{E,g,b}$. 
Due to ellipticity of $\Delta_{E,g,b}$, the space $\Omega_{g,b}(M;E)(0)$ is finite dimensional and consists of smooth forms only.
Since $d_E$ commutes with $\Delta_{E,g,b}$ the zero eigen space $\Omega_{g,b}(M;E)(0)$ is a subcomplex of $\Omega(M;E)$.
As $\Delta_{E,g,b}$ fails to be selfadjoint,
the differential on $\Omega_{g,b}(M;E)(0)$ will in general not vanish. However, the
inclusion $\Omega_{g,b}(M;E)(0)\to\Omega(M;E)$ induces an isomorphism in cohomology. It follows that the restriction
of \eqref{E:Int} to $\Omega_{g,b}(M;E)(0)$ induces an isomorphism in cohomology too.
Because $\Delta_{E,g,b}$ is symmetric with respect to $\beta$, the bilinear form $\beta$ will 
restrict to a non-degenerate symmetric bilinear form on $\Omega_{g,b}(M;E)(0)$. 
We thus have a quasi isomorphism between finite dimensional complexes
$$
\Int_{X,E}:\Omega_{g,b}(M;E)(0)\to C_b(X;E)
$$
each of which is equipped with a non-degenerate symmetric bilinear form. This permits to define the square of the relative torsion which we will
denote by
$$
\tau\Bigl(\Omega_{g,b}(M;E)(0)\xrightarrow{\Int_{X,E}}C_b(X;E)\Bigr)\in\C^\times.
$$
A more detailed discussion of these facts can be found in \cite[Section~4]{BH06b}.

Recall the Mathai--Quillen form, see \cite{MQ86} or \cite[Section~III]{BZ92},
$$
\Psi_g\in\Omega^{n-1}(TM\setminus M;\mathcal O_M),
$$ 
and the Kamber--Tondeur form
\begin{equation}\label{E:KTT} 
\omega_{E,b}:=-\tfrac12\tr_E(b^{-1}\nabla^Eb)\in\Omega^1(M;\C).
\end{equation} 
Note that $\omega_{E,b}$ is closed since $\nabla^E$ is flat, see \cite[Section~2]{BH06b}. The integral
$$
\int_{M\setminus\mathcal X}\omega_{E,b}\wedge(-X)^*\Psi_g
$$
will in general not converge, but can easily be regularized, see \cite[Section~2]{BH03}, \cite[Section~3]{BH06} or \cite[Section~III]{BZ92}.
Consider the non-vanishing complex number
\begin{multline*}
\mathcal S_{E,g,b,X}:=\tau\Bigl(\Omega_{g,b}(M;E)(0)\xrightarrow{\Int_{X,E}} C_b(X;E)\Bigr)
\\\cdot\prod_q\bigl({\det}'(\Delta_{E,g,b,q})\bigr)^{(-1)^qq}
\cdot\exp\Bigl(-2\int_{M\setminus\mathcal X}\omega_{E,g}\wedge(-X)^*\Psi_g\Bigr).
\end{multline*}
In \cite{BH06b} the following result, analogous to the anomaly formula for the classical Ray--Singer torsion \cite[Theorem~0.1]{BZ92}, 
has been established.\footnote{Strictly speaking, this was done for vanishing Euler--Poincar\'e characteristics only. However,
with few additional elementary arguments Theorem~\ref{T:anom} below can be proved exactly as in \cite[Section~6]{BH06b}, the
crucial analytic results \cite[Proposition~6.1 and 6.2]{BH06b} have been established without any restriction on the 
Euler--Poincar\'e characteristics.}

\begin{theorem}\label{T:anom}
The quantity $\mathcal S_{E,g,b,X}$ is independent of the Morse--Smale vector field $X$, independent of the
Riemannian metric $g$ and locally constant in $b$. It thus depends on the flat bundle $E$ and the homotopy class $[b]$
of the fiber wise non-degenerate bilinear form only, and will be denoted by $\mathcal S_{E,[b]}$.
\end{theorem}

\begin{remark}\label{R:Scomban}
There is a conceptual interpretation of $\mathcal S_{E,[b]}$ as the quotient of two invariants, see \cite[Section~5]{BH06b}.
\end{remark}

In analogy with a result of Cheeger \cite{C77,C79}, M\"uller \cite{M78} and Bismut--Zhang \cite[Theorem~0.2]{BZ92}, 
the following conjecture was raised in \cite[Conjecture~5.1]{BH06b}.

\begin{conjecture}\label{C:main}
We have $\mathcal S_{E,[b]}=1$ for every flat complex vector bundle $E$ and every fiber wise non-degenerate symmetric bilinear form $b$ on $E$.
\end{conjecture}

This conjecture has been verified in several non-trivial situations, see \cite{BK06b} and  \cite[Section~5]{BH06b}.
One purpose of this paper is to establish Conjecture~\ref{C:main} for odd dimensional manifolds, up to sign.
More precisely, we will show

\begin{theorem}\label{T:main}
Suppose $M$ is odd dimensional. Then $\mathcal S_{E,[b]}=\pm1$ for
every flat complex vector bundle $E$ and every fiber wise 
non-degenerate symmetric bilinear form $b$ on $E$.\footnote{In the appendix we show how one can remove the
sign ambiguity by extending Witten--Helffer--Sj\"ostrand theory to generalized Morse functions, a project 
partially realized in \cite{HKW95}. One also note that an extension of the theorem above to compact manifolds with boundary  implies the result for closed  even dimensional manifolds}
\end{theorem}

In fact we actually show that the strategy of proving the Cheeger--M\"uller theorem presented in \cite{BFK96} works here too,
almost identically. However, the spectral properties of the Witten deformation of the Laplace operators associated with a 
Riemannian metric and a non-degenerate symmetric bilinear form can not benefit from the methods described in more details in \cite {BFKM96},
based on selfadjointness.
Fortunately, the key geometric consequences continue to hold. The other purpose of this paper is to extend
Witten--Helffer--Sj\"ostrand theory to this non-selfadjoint situation, which may be of independent interest.

This paper is organized as follows.
In section~\ref{S:whs} we will extend Witten--Helffer--Sj\"ostrand theory \cite{BZ92,BZ94, BFK96,BFK01,BFKM96,HS84,HS85} to this non-selfadjoint situation.
More precisely, we choose a Morse function $f$ on $M$ and fix a Riemannian metric $g$ which has a standard form near the critical points $\mathcal X$ of $f$.
We assume that the gradient vector field $X:=-\grad_g(f)$ satisfies the Smale transversality condition. Finally, we assume that the bilinear
form $b$ is parallel in a neighbourhood of $\mathcal X$, with respect to $\nabla^E$. In view of Theorem~\ref{T:anom} these assumptions
do not cause a loss of generality for the purpose of computing $\mathcal S_{E,[b]}$. We then consider the family of Witten deformed flat connections on $E$
\begin{equation}\label{E:XX2}
\nabla^{E_u}:=\nabla^E+udf,\qquad u\geq0.
\end{equation}
Let us write $E_u$ for the complex vector bundle $E$ equipped with the flat connection $\nabla^E+udf$. 
Since $e^{uf}:(E_u,b)\to(E,e^{-2uf}b)$ is an isomorphism of flat vector bundles with bilinear forms, it follows from Theorem~\ref{T:anom} that
$\mathcal S_{E_u,[b]}$ is constant in $u$.

Let us introduce the \emph{large analytic torsion}
\begin{equation}\label{E:laanto}
\tau_{\la,u}:=\prod_q\bigl({\det}^\la(\Delta_{E_u,g,b,q})\bigr)^{(-1)^qq}
\end{equation}
where ${\det}^\la(\Delta_{E,g,b,q})$ denotes the zeta regularized product of eigen values whose real part is larger than $1$.
In view of Proposition~\ref{P:specgap} below, $\tau_{\la,u}$ is analytic in $u$ for sufficiently large $u$. 
Moreover, let us write $\Omega_\sm(M;E_u)$ for the sum of eigen spaces corresponding to eigen values with real part at most $1$.
We refer to $\Omega_\sm(M;E_u)$ as the \emph{small complex,} see Proposition~\ref{P:specgap} below.
This is a finite dimensional subcomplex of $\Omega(M;E_u)$, $\beta$ restricts to a non-degenerate symmetric bilinear form $\beta_{\sm,u}$
on $\Omega_\sm(M;E_u)$, and the restriction of \eqref{E:Int} to $\Omega_\sm(M;E_u)$ provides a quasi isomorphism
\begin{equation}\label{E:intusm1}
\Int_{\sm,u}\co\Omega_\sm(M;E_u)\to C(X;E_u).
\end{equation}
Let us write $\tau(\Int_{\sm,u})\in\C^\times$ for the relative torsion of \eqref{E:intusm1}.
It is not hard to show, see \cite[Proposition~5.10]{BH06b}, that
\begin{equation}\label{E:DD}
\mathcal S_{E,[b]}=\mathcal S_{E_u,[b]}=\tau(\Int_{\sm,u})\cdot\tau_{\la,u}\cdot\exp\Bigl(-2\int_{M\setminus\mathcal X}\omega_{E_u,b}\wedge(-X)^*\Psi_g\Bigr).
\end{equation} 
The Witten--Helffer--Sj\"ostrand estimates show that for sufficiently large $u$ the integration \eqref{E:intusm1} is an isomorphism
of complexes, and they provide an asymptotic comparison of the bilinear form $\beta_{\sm,u}$ on $\Omega_\sm(M;E_u)$ and 
the bilinear form $b_{\mathcal X}$ on $C(X;E_u)$. The precise statement is contained in Theorem~\ref{T:WHS} below. This result permits to
compute the asymptotic expansion of $\tau(\Int_{\sm,u})$ as $u\to\infty$, see Corollary~\ref{C:WHS}. 
In view of \eqref{E:DD} this yields a formula for $\mathcal S_{E,[b]}$
in terms of the free  term of the asymptotic expansion of $\tau_{\la,u}$ as $u\to\infty$, see Corollary~\ref{C:WHSla}.

Unfortunately we are unable at this time to calculate directly this constant term and check whether $\mathcal S_{E,[b]}$ is one or not.
However, in Section~\ref{S:assla} we show that for two systems $(M,E,g,b,f)$ and $(\tilde M,\tilde E,\tilde g,\tilde b,\tilde f)$
with $M$ and $\tilde M$ of the same dimension, $E$ and $\tilde E$ of the same rank, $f$ and $\tilde f$ with the same number of critical points
in each index
$$
\log\tau_{\la,u}-\log\tilde\tau_{\la,u}
$$
has an asymptotic expansion whose free term is computable as integral of local quantities, see Theorem~\ref{T:rella}.
This is done as in \cite{BFK96}; precisely, combining the fact that the Witten deformation is elliptic with parameter
away from the critical points with a Mayer--Vietoris formula for the zeta regularized determinants of elliptic operators,
yields a result as formulated in Theorem~\ref{T:rella}.
Note that we have an unambiguously defined logarithm, given by the formula:
$$
\log\tau_{\la,u}:=\sum_q(-1)^qq\frac\partial{\partial s}\Big|_{s=0}\sum_{\lambda\in\Spec\Delta_{u,q},\Re\lambda>1}\lambda^{-s}
$$
Playing with its symmetry, as in \cite{BFK96}, we derive that, for odd dimensional manifolds,
$$
\mathcal S_{E,[b]}^2=\mathcal S_{\tilde E,[\tilde b]}^2.
$$
We will then use this to give a proof of Theorem~\ref{T:main}.

\section{Witten--Helffer--Sj\"ostrand theory}\label{S:whs}

Let $f$ be a Morse function on a closed connected smooth manifold $M$ of dimension $n$.
Fix a Riemannian metric $g$ on $M$
so that the vector field $X:=-\grad_g(f)$ satisfies the
Smale transversality condition. Let $\mathcal X\subseteq M$ denote the set of 
critical points of $f$, and denote by $\ind(x)\in\{0,1,\dotsc,n\}$ 
the Morse index of a critical point $x\in\mathcal X$. For every critical point
$x\in\mathcal X$ we fix an open neighbourhood $U_x$ of $x$
and a diffeomorphisms (Morse chart)
$\varphi_x=(\varphi_x^1,\dotsc,\varphi_x^n)\co U_x\to\R^n$ so that $\varphi_x(x)=0$
and
\begin{equation}\label{E:M1}
f=f(x)-\frac12\sum_{i\leq\ind(x)}(\varphi_x^i)^2+\frac12\sum_{i>\ind(x)}(\varphi_x^i)^2
\end{equation}
in a neighbourhood of $x$. 
We will assume that every $x\in\mathcal X$ admits a neighbourhood
on which the Riemannian metric takes the form
\begin{equation}\label{E:M2}
g=\sum_{i=1}^nd\varphi_x^i\otimes d\varphi_x^i.
\end{equation}
Finally, we will assume that, in a neighbourhood of $\mathcal X$, we have 
\begin{equation}\label{E:M3}
\nabla^Eb=0.
\end{equation}

Recall the family of integration homomorphisms \eqref{E:intusm1} associated with the Witten deformed
flat bundles $E_u$, see \eqref{E:XX2}. For $u>0$ let us introduce the scaling isomorphism 
\begin{equation}\label{E:scaling}
\eta_u\co C(X;E_u)\to C(X;E_u)
\end{equation}
defined by
$$
\eta_u(w):=\Bigl(\frac\pi u\Bigr)^{n/4-q/2}w,\qquad w\in C^q(X;E_u).
$$
The aim of this section is to provide a proof of the following

\begin{theorem}\label{T:WHS}
For sufficiently large $u$, the restriction of the integration map
\begin{equation}\label{E:intusm}
\Int_{u,\sm}\co\Omega_\sm(M;E_u)\to C(X;E_u)
\end{equation}
is an isomorphism of complexes. Moreover, there exists a constant $\epsilon>0$ 
so that
$$
(\eta_u\Int_{u,\sm})_*\beta_{u,\sm}=b_{\mathcal X}+\LO(e^{-\epsilon u})
\qquad\text{as $u\to\infty$.}
$$
\end{theorem}

Theorem~\ref{T:WHS} generalizes the classical Witten--Helffer--Sj\"ostrand theorem for selfadjoint Laplacians, see
\cite{HS84}, \cite{HS85}, \cite{BZ92}, \cite{BZ94} and \cite[Theorem~5.5]{BFKM96},
to this non-selfadjoint situation. The proof of this result will be similar to the one in \cite{BFKM96}. 
A few additional arguments, however, are necessary since the
Laplacians $\Delta_u:=\Delta_{E_u,g,b}$ are not necessarily selfadjoint. A key step
is to establish a widening gap in the spectrum of $\Delta_u$, as $u\to\infty$.
A finite number of eigen values will exponentially fast approach $0$,
while the real part of the remaining will grow linearly with $u$. For
the precise statement see Proposition~\ref{P:specgap} below.
This justifies the term \emph{small complex.}

In order to spell out two corollaries let us introduce the notation
\begin{align*}
\chi&:=\sum_q(-1)^q\dim C^q(X;E)=\chi(M)\rank(E)
\\
\chi'&:=\sum_q(-1)^qq\dim C^q(X;E)=\sum_q(-1)^qq|\mathcal X_q|\rank(E) 
\end{align*}
where $|\mathcal X_q|$ denotes the number of critical points of index $q$.

\begin{corollary}\label{C:WHS}
There exists a constant $\epsilon>0$ so that, as $u\to\infty$,
$$
\tau(\Int_{u,\sm})=\Bigl(\frac\pi u\Bigr)^{\frac n2\chi-\chi'}\bigl(1+\LO(e^{-\epsilon u})\bigr).
$$
\end{corollary}

\begin{proof}
In view of Theorem~\ref{T:WHS} the integration \eqref{E:intusm} is an 
isomorphism of complexes, assuming $u$ is sufficiently large. Hence, we have
\begin{equation}\label{E:AA}
\tau(\Int_{u,\sm})=\sdet\Bigl((b_{\mathcal X})^{-1}\bigl((\Int_{u,\sm})_*\beta_{u,\sm}\bigr)\Bigr).
\end{equation}
Here the right hand side denotes the super determinant\footnote{
If $\varphi_k:V^k\to V^k$ is a linear endomorphism of a graded vector space which preserves the grading, then its
super (or graded) determinant is given by $\sdet(\varphi)=\prod_k(\det\varphi_k)^{(-1)^k}$.}  
of the composition
$$
C(X;E_u)\xrightarrow{(\Int_{u,\sm})_*\beta_{u,\sm}}C(X;E_u)'\xrightarrow{(b_{\mathcal X})^{-1}}C(X;E_u)
$$
and the non-degenerate bilinear forms are considered as isomorphisms between
the vector space $C(X;E_u)$ and its dual, $C(X;E_u)'$.
Note that as isomorphisms $C(X;E_u)\to C(X;E_u)'$ we have
$$
(\Int_{u,\sm})_*\beta_{u,\sm}=\bigl((\eta_u\Int_{u,\sm})_*\beta_{u,\sm}\bigr)\circ\eta_u^2
$$
From \eqref{E:AA} we thus conclude
\begin{equation}\label{E:BB}
\tau(\Int_{u,\sm})
=\sdet\Bigl((b_{\mathcal X})^{-1}\bigl((\eta_u\Int_{u,\sm})_*\beta_{u,\sm}\bigr)\Bigr)
\cdot(\sdet\eta_u)^2.
\end{equation}
From Theorem~\ref{T:WHS} we obtain a constant $\epsilon>0$ so that, as $u\to\infty$,
$$
\sdet\Bigl((b_{\mathcal X})^{-1}\bigl((\eta_u\Int_{u,\sm})_*\beta_{u,\sm}\bigr)\Bigr)
=1+\LO(e^{-\epsilon u}).
$$
One readily checks:
$$
(\sdet\eta_u)^2=\Bigl(\frac\pi u\Bigr)^{\frac n2\chi-\chi'}
$$
Combining the latter two with \eqref{E:BB} completes the proof of the corollary.
\end{proof}

\begin{corollary}\label{C:WHSla}
There exist (unique) constants $a_0\in\C^\times$, $a_1,a_2\in\C$ with the following property.
There exists a constant $\epsilon>0$ so that, as $u\to\infty$,
$$ 
\tau_{\la,u}=a_0\cdot e^{a_1u}\cdot u^{a_2}\cdot\bigl(1+\LO(e^{-\epsilon u})\bigr).
$$
These constants are given by:
\begin{align*}
a_0&=\mathcal S_{E,[b]}\cdot\pi^{-(\frac n2\chi-\chi')}\cdot\exp\Bigl(2\int_{M\setminus\mathcal X}\omega_{E,b}\wedge(-X)^*\Psi_g\Bigr)
\\
a_1&=2\rank(E)\int_{M\setminus\mathcal X}df\wedge(-X)^*\Psi_g
\\
a_2&=\tfrac n2\chi-\chi'
\end{align*}
\end{corollary}

\begin{proof}
Combining Corollary~\ref{C:WHS} with \eqref{E:DD} we obtain a constant $\epsilon>0$ so that
\begin{equation}\label{E:1010}
\tau_{\la,u}=\mathcal S_{E,[b]}\cdot
\Bigl(\frac u\pi\Bigr)^{\frac n2\chi-\chi'}
\cdot\exp\Bigl(2\int_{M\setminus\mathcal X}\omega_{E_u,b}\wedge(-X)^*\Psi_g\Bigr)\cdot
\bigl(1+\LO(e^{-\epsilon u})\bigr)
\end{equation}
as $u\to\infty$. From $\nabla^{E_u}b=\nabla^Eb-2udf\otimes b$, see \eqref{E:XX2}, we obtain
\begin{equation}\label{E:KTpre}
b^{-1}\nabla^{E_u}b=b^{-1}\nabla^Eb-2udf\otimes1_E
\end{equation}
Therefore, see \eqref{E:KTT}, $\omega_{E_u,b}=\omega_{E,b}+u\rank(E)df$ and thus
\begin{multline*}
\int_{M\setminus\mathcal  X}\omega_{E_u,b}\wedge(-X)^*\Psi_g
\\=\int_{M\setminus\mathcal X}\omega_{E,b}\wedge(-X)^*\Psi_g
+u\rank(E)\int_{M\setminus\mathcal X}df\wedge(-X)^*\Psi_g.
\end{multline*}
Combining this with \eqref{E:1010} the corollary follows.
\end{proof}

\subsection*{A formula for the Witten perturbed Laplacian}

Recall that for every critical point $x\in\mathcal X$ we have fixed a Morse chart
$$
\varphi_x=(\varphi_x^1,\dotsc,\varphi_x^n)\co U_x\to\R^n.
$$ 
Choose $\rho>0$ so that with 
$$
B_x:=\varphi_x^{-1}\bigl(\{z\in\R^n\mid|z|<\rho\}\bigr),\qquad x\in\mathcal X,
$$
equations \eqref{E:M1}, \eqref{E:M2} and \eqref{E:M3} hold on $B_x$, for every
critical point $x\in\mathcal X$. We assume $\rho$ was chosen sufficiently small so
that the closures of $B_x$ are mutually disjoint. 
It will be convenient to further assume, by choosing $\rho$ sufficiently
small, that for all $x,y\in\mathcal X$ with $\ind(x)=\ind(y)$, we have
\begin{equation}\label{E:SmB}
W_x^-\cap B_y=
\begin{cases}
\emptyset&\text{if $x\neq y$}
\\\varphi_x^{-1}(\R^{\ind(x)})\cap B_x&\text{if $x=y$}
\end{cases}
\end{equation}
This is possible in view of the Smale transversality property of $X$.
Here $W_x^-$ denotes the unstable manifold of $x$.
Set
\begin{equation}\label{E:B}
B:=\bigcup_{x\in\mathcal X}B_x.
\end{equation}
For each $x\in\mathcal X$ we introduce the smooth radial function
\begin{equation}\label{E:rx2}
r_x^2\co B_x\to\R,\qquad r_x^2:=\sum_{i=1}^n(\varphi_x^i)^2.
\end{equation}

Over $B_x$, $x\in\mathcal X$, we decompose the cotangent bundle as 
$$
T^*M|_{B_x}=V_x^-\oplus V_x^+,
$$
where the subbundle $V_x^-\subseteq T^*M|_{B_x}$ is spanned by 
$d\varphi_x^i$, $1\leq i\leq\ind(x)$, and the subbundle $V_x^+\subseteq T^*M|_{B_x}$ 
is spanned by $d\varphi_x^i$, $\ind(x)<i\leq n$. Let us write $\Lambda:=\Lambda^*T^*M$.
We obtain an induced decomposition 
$$
\Lambda|_{B_x}=\Lambda_x^-\otimes\Lambda_x^+,\qquad\text{where}\qquad
\Lambda_x^\pm:=\Lambda^*V_x^\pm.
$$
Let us write 
$$
N^\pm_x\in\Gamma(\eend(\Lambda_x^\pm))
$$ 
for the grading operator acting by multiplication with $q$ on
$\Lambda^qV_x^\pm$. We will denote the operators
$N_x^-\otimes\id_{\Lambda_x^+}\otimes\id_E$ and $\id_{\Lambda_x^-}\otimes N^+_x\otimes\id_E$ 
acting on $\Lambda\otimes E|_{B_x}$ by $N_x^-$ and $N_x^+$ too.

\begin{lemma}\label{L:F1}
There exists a zero order operator $L\in\Gamma(\eend(\Lambda\otimes E))$
so that for all $u\geq0$ we have
$$
\Delta_u=\Delta_0+uL+u^2|df|^2.
$$
For every critical point $x\in\mathcal X$ we have
$$
L|_{B_x}=2\bigl(N^+_x+\ind(x)-N^-_x\bigr)-n,
$$
and, for $u\geq0$,
$$
\Delta_u=\Delta_0-un+u^2r_x^2+2u\bigl(N^+_x+\ind(x)-N^-_x\bigr)\qquad\text{over $B_x$.}
$$
\end{lemma}

\begin{proof}
As in the selfadjoint situation this can be verified in coordinates, see \cite{HS85}.
Alternatively, one can give a conceptual proof based on the observation, see \cite[Section~IV]{BZ92}, that the Laplacians
$\Delta_u$ are the squares of Dirac operators associated to Clifford super connections on $\Lambda\otimes E$.
\end{proof}

\begin{remark}\label{R:F2}
Let $x\in\mathcal X$ be a critical point. Use the flat connection $\nabla^E$
and the (flat) Levi--Civita connection to identify 
$\Omega(B_x;E)=C^\infty(B_x,\Lambda\otimes E_x)$. Then, via this 
identification, we have
$$
\Delta_0v=-\sum_i\frac{\partial^2v}{(\partial\varphi_x^i)^2},
\qquad v\in C^\infty(B_x,\Lambda\otimes E_x).
$$
\end{remark}

\subsection*{Compatible Hermitian structure}

We will now introduce a Hermitian structure on the vector bundle $E$.
It is possible to choose such a Hermitian structure to be compatible with the
symmetric bilinear form. To this end we start with

\begin{lemma}\label{LLL:1}
There exists a fiber wise complex anti-linear involution $v\mapsto\bar v$ on
$E$ such that, for $e_1$ and $e_2$ in the same fiber of $E$, $e_1\neq0$, we have
$$
b(\bar e_1,\bar e_2)=\overline{b(e_1,e_2)}\qquad\text{and}\qquad b(e_1,\bar e_1)>0.
$$
Moreover, this involution can be chosen to be parallel over $B$, 
see \eqref{E:B}. That is, for $e\in E|_B$ we have
$$
\nabla^E\bar e-\overline{\nabla^Ee}=0.
$$ 
\end{lemma}

\begin{proof}
The fiber wise non-degenerate symmetric bilinear form $b$ provides a reduction
of the structure group to $O_k(\C)$, where $k=\rank(E)$.
The natural inclusion $O_k(\R)\to O_k(\C)$ is a homotopy equivalence,
hence the structure group can be further reduced to $O_k(\R)$.
Note that the subgroup $O_k(\R)\subseteq O_k(\C)$ consists of those
matrices whose action on $\C^k$ commutes with the standard complex
conjugation on $\C^k$. The existence of the desired complex 
anti-linear involution on $E$ follows immediately.
\end{proof}

We fix a complex conjugation as in Lemma~\ref{LLL:1}. Then
\begin{equation}\label{EEEE:1}
\langle e_1,e_2\rangle:=b(e_1,\bar e_2)
\end{equation}
defines a fiber wise (positive definite) Hermitian structure on $E$.
We will write $|e|$ for the associated fiber wise norm.
Note that 
\begin{equation}\label{EEEE:1.5}
b(\bar e_1,\bar e_2)=\overline{b(e_1,e_2)},\quad
\langle\bar e_1,\bar e_2\rangle=\overline{\langle e_1,e_2\rangle}
\quad\text{and}\quad|\bar e|=|e|.
\end{equation}
Moreover, this Hermitian structure is parallel over $B$, with respect to $\nabla^E$.

\begin{remark}
Note that $F:=\{e\in E\mid \bar e=e\}$ is a real subbundle of $E$,
and that there is a canonical isomorphism $F\otimes\C=E$. The restrictions of $b$
and $\langle\cdot,\cdot\rangle$ to $F$ coincide, and define a (positive
definite) Euclidean structure on $F$. We can understand both, $b$ and
$\langle\cdot,\cdot\rangle$, as complexifications of this Euclidean 
inner product, once complexifying to a bilinear form and once complexifying
to a sesquilinear form.
\end{remark}

The fiber wise Hermitian structure on $E$ induces a Hermitian inner product on the
Morse complex $C(X;E_u)$ in an obvious way. For $a_1,a_2\in C(X;E_u)$
we will denote this by $\langle a_1,a_2\rangle_{\mathcal X}$. 
Similarly we will write $|a|_{\mathcal X}$ for the
associated norm, $a\in C(X;E_u)$. Moreover, the fiber wise complex conjugation on $E$
induces a complex conjugation on $C(X;E_u)$. For $a,a_1,a_2\in C(X;E_u)$
we have, see \eqref{EEEE:1} and \eqref{EEEE:1.5},
\begin{equation}\label{E:A50}
\langle a_1,a_2\rangle_{\mathcal X}=b_{\mathcal X}(a_1,\bar a_2)
\end{equation}
as well as
\begin{equation}\label{E:A51}
b_{\mathcal X}(\bar a_1,\bar a_2)=\overline{b_{\mathcal X}(a_1,a_2)},
\quad
\langle\bar a_1,\bar a_2\rangle_{\mathcal X}=\overline{\langle a_1,a_2\rangle_{\mathcal X}}
\quad\text{and}\quad
|\bar a|_{\mathcal X}=|a|_{\mathcal X}.
\end{equation}

Using the Riemannian metric $g$, we obtain an induced fiber wise Hermitian inner
product on $\Lambda^*T^*M\otimes E$ which will be denoted by 
$\langle v,w\rangle_g$, $v,w\in\Omega(M;E_u)$. Then
$$
\llangle v,w\rrangle:=\int_M\langle v,w\rangle_g\vvol_g,\qquad v,w\in\Omega(M;E_u)
$$
is a Hermitian inner product on $\Omega(M;E_u)$. We will write $\|v\|$ for the associated $L_2$--norm, $v\in\Omega(M;E_u)$.
The complex conjugation induces a complex conjugation on $\Omega(M;E_u)$.
For $v,w\in\Omega(M;E_u)$ we clearly have, see \eqref{EEEE:1} and \eqref{EEEE:1.5},
\begin{equation}\label{EEEE:2} 
\llangle v,w\rrangle=\beta(v,\bar w)
\end{equation}
as well as
\begin{equation}\label{EEEE:4}
\beta(\bar v,\bar w)=\overline{\beta(v,w)},\quad
\llangle\bar v,\bar w\rrangle=\overline{\llangle v,w\rrangle}\quad\text{and}\quad
\|\bar v\|=\|v\|.
\end{equation}
From $\beta(\Delta_u v,w)=\beta(v,\Delta_u w)$ we thus also obtain
\begin{equation}\label{EEEE:3}
\llangle\Delta_u v,w\rrangle
=\llangle\Delta_u\bar w,\bar v\rrangle
=\llangle v,\overline{\Delta_u\bar w}\rrangle,\qquad v,w\in\Omega(M;E_u).
\end{equation}

\begin{remark}\label{R:adjoint}
It follows from Lemma~\ref{L:F1} that, over $B$, the Laplacian $\Delta_u$ commutes with the complex
conjugation, that is
$$
\Delta_u\bar w=\overline{\Delta_uw},\qquad w\in\Omega(M;E_u),\ \supp w\subseteq B.
$$
From \eqref{EEEE:3} we thus get
$$
\llangle\Delta_uv,w\rrangle=\llangle v,\Delta_uw\rrangle,\qquad v,w\in\Omega(M;E_u),\ \supp w\subseteq B.
$$
It also follows from Lemma~\ref{L:F1} that over $B$, the Laplacian $\Delta_u$ coincides with the selfadjoint Witten Laplacian
associated to the compatible Hermitian structure, see \cite{BZ92} or \cite{BFKM96}.
\end{remark}

\subsection*{Construction of approximate small eigen forms}

In this section we will construct as in \cite {HS84}, cf.\ also \cite{BFKM96}, finite dimensional subspaces
$V_u\subseteq\Omega(M;E_u)$ which approximate $\Omega_\sm(M;E_u)$ as $u\to\infty$.
The estimates in the section show that $\Delta_u$ has small norm on $V_u$
and is large on the orthogonal complement of $V_u$. 

For a critical point $x\in\mathcal X$ and $e\in E_x$ we let $\tilde e\in\Gamma(E|_{B_x})$
denote the unique parallel section, i.e.\ $\nabla^E\tilde e=0$, satisfying 
$\tilde e(x)=e$. Moreover, we introduce the differential form
$$
\Omega_x^-:=d\varphi_x^1\wedge\cdots\wedge d\varphi_x^{\ind(x)}\in\Omega^{\ind(x)}(B_x;\R).
$$
Choose a smooth function $\sigma\co\R\to[0,1]$ such that
$\sigma(t)=1$ for all $t\leq \rho/3$ and $\sigma(t)=0$ for all $t\geq2\rho/3$.
For $u\geq0$, consider the smooth form, see \eqref{E:rx2},
\begin{equation}\label{E:phiuxe}
\phi_{u,e}:=(\sigma\circ r_x)e^{-ur_x^2/2}\Omega^-_x\otimes\tilde e\in\Omega(B_x;E_u).
\end{equation}  
Since $\phi_{u,e}$ has compact support contained in $B_x$, 
we can consider it as a globally defined form, $\phi_{u,e}\in\Omega(M;E_u)$.

\begin{definition}\label{D:Vu}
For $u\geq0$ we let $V_u\subseteq\Omega(M;E_u)$ denote the finite dimensional
subspace spanned by the forms $\phi_{u,e}$ where $x$ runs through $\mathcal X$
and $e$ runs through (a basis) of $E_x$.
\end{definition}

\begin{observation}\label{O:Vubar}
The complex conjugation preserves $V_u$. The $\beta$- orthogonal complement
of $V_u$ coincides with its Hermitian orthogonal complement. It will be denoted
by $V_u^\perp\subseteq\Omega(M;E_u)$.
\end{observation}

\begin{proof}
The first assertion is immediate from the definition of $V_u$ and the 
fact that the complex conjugation is parallel over $B$, see Lemma~\ref{LLL:1}. 
The second claim then follows from \eqref{EEEE:2}. 
\end{proof}

For $k\in\N$ let $\|w\|_{C^k}$ denote the (a fixed) $C^k$--norm of $w\in\Omega(M;E_u)$. 
The following estimates follow easily from the  structure of $\Delta_u$ in the neighborhood of critical 
points, cf.\ Lemma~\ref{L:F1} and Remark~\ref{R:F2}.

\begin{lemma}\label{L:AAA}
There exist constants $u_{1,k}\geq1$ and $\epsilon_{1,k}>0$ so that for all
$k\in\N$, $u\geq u_{1,k}$ and $v\in V_u$ we have
$$
\|\Delta_uv\|_{C^k}\leq e^{-\epsilon_{1,k}u}\|v\|.
$$
\end{lemma}

\begin{proof}
Let $x\in\mathcal X$ and $e\in E_x$. An elementary computation, see Remark~\ref{R:F2}, 
shows
$$
(\Delta_0-un+u^2r_x^2)e^{-ur_x^2/2}\Omega^-_x\otimes\tilde e=0.
$$
We conclude, see \eqref{E:phiuxe}, that 
$$
\supp\bigl((\Delta_0-un+u^2r_x^2)\phi_{u,e}\bigr)\subseteq R_x,
$$ 
where $R_x:=\varphi_x^{-1}\bigl(\{z\in\R^n\mid \rho/3\leq|z|\leq2\rho/3\}\bigr)$.
Note that there exist constants $C_k\geq0$ and $\epsilon'>0$ 
so that, for all $u\geq0$,
$$
\bigl\|(\Delta_0-un+u^2r_x^2)\phi_{u,e}\bigr\|_{C^k}\leq C_ke^{-\epsilon'u}|e|.
$$
Together with $(N^+_x+\ind(x)-N_x^-)\Omega_x^-=0$ and the formula 
in Lemma~\ref{L:F1} we see that, for all $u\geq0$,
\begin{equation}\label{E:Y4}
\bigl\|\Delta_u\phi_{u,e}\bigr\|_{C^k}\leq C_ke^{-\epsilon'u}|e|.
\end{equation}
Clearly,
$$
\|\phi_{u,e}\|^2=\int_{B_x}(\sigma\circ r_x)^2e^{-ur_x^2}|e|^2\vvol_g
=\int_{\R^n}\sigma(|z|)^2e^{-u|z|^2}dz|e|^2,
$$
and so there exist constants $C\geq0$ and $\epsilon''>0$ so that for all $u\geq1$
$$
\Bigl|\|\phi_{u,e}\|-(\pi/u)^{n/4}|e|\Bigr|\leq Ce^{-\epsilon''u}|e|.
$$
Combining this with \eqref{E:Y4} we find constants $u_k\geq1$ and $\epsilon_{1,k}>0$
so that for $u\geq u_k$, $x\in\mathcal X$ and $e\in E_x$ we have
$$
\bigl\|\Delta_u\phi_{u,e}\|_{C^k}\leq e^{-\epsilon_{1,k}u}\|\phi_{u,e}\|.
$$
The lemma now follows  from the fact that for $x_1,x_2\in\mathcal X$, $x_1\neq x_2$,
and $e_1\in E_{x_1}$, $e_2\in E_{x_2}$ the forms $\phi_{u,e_1}$ and $\phi_{u,e_2}$
have disjoint support. 
\end{proof}

\begin{lemma}\label{L:BBB}
There exist constants $u_2\geq1$ and $\epsilon_2>0$ so that for $u\geq u_2$,
$v\in V_u$ and $w\in\Omega(M;E_u)$ we have
$$
|\llangle\Delta_uv,w\rrangle|\leq e^{-\epsilon_2u}\|v\|\|w\|
\qquad\text{and}\qquad
|\llangle v,\Delta_uw\rrangle|\leq e^{-\epsilon_2u}\|v\|\|w\|.
$$
\end{lemma}

\begin{proof}
The first statement follows immediately from the Cauchy--Schwarz inequality
and Lemma~\ref{L:AAA} for $k=0$. To see the second inequality, recall from
Observation~\ref{O:Vubar} that $V_u$ is invariant under the complex conjugation.
Hence \eqref{EEEE:3}, the first statement and \eqref{EEEE:4} imply
\begin{equation*}
|\llangle v,\Delta_uw\rrangle|=|\llangle\Delta_u\bar v,\bar w\rrangle|
\leq e^{-\epsilon_2u}\|\bar v\|\|\bar w\|
=e^{-\epsilon_2u}\|v\|\|w\|.
\qedhere
\end{equation*}
\end{proof}

Introduce the smooth cut-off function, see \eqref{E:rx2} and \eqref{E:phiuxe},
$$
\chi\co M\to[0,1],\quad
\chi(y):=\sum_{x\in\mathcal X}\sigma\circ r_x.
$$
Note that $\chi=1$ in a neighbourhood of $\mathcal X$ and $\supp\chi\subseteq B$, see \eqref{E:B}.

\begin{lemma}\label{L:2}
There exist constants $u_3'\geq1$ and $\epsilon_3'>0$ so that for any $u\geq u_3'$ and $v'\in V^\perp_u$ we have
\begin{equation}\label{E:101}
\Re\llangle\Delta_u(\chi v'),\chi v'\rrangle \geq \epsilon_3'u \|\chi v'\|^2.
\end{equation}
\end{lemma}

\begin{proof}
It suffices to establish the result  in $\R^n$ and for the standard Witten Laplacian only, see
Remark~\ref{R:adjoint}. This is done in \cite[Appendix~A.2, equation~(5.7)]{BFKM96}.
\end{proof}

\begin{lemma}\label{L:4}
There exist constants $u_3''\geq1$ and $\epsilon_3''>0$ so that for all $w\in\Omega(M;E_u)$ with
$\supp w\subseteq\overline{\{y\in M\mid\chi(y)\neq1\}}$ we have
$$
\Re\llangle\Delta_uw,w\rrangle\geq\epsilon_3''u^2\|w\|^2.
$$
\end{lemma}

\begin{proof}
Since $\Delta_0+\Delta_0^*$ is a selfadjoint operator whose principal symbol is positive definite
it follows, see \cite[Corollary~9.3]{S01}, that $\Delta_0+\Delta_0^*$ is bounded from below, i.e.\
there exists a constant $C\geq0$ so that for all $w\in\Omega(M;E_u)$
\begin{equation}\label{E:B12}
2\Re\llangle\Delta_0w,w\rrangle=\llangle(\Delta_0+\Delta_0^*)w,w\rrangle\geq-C\|w\|^2.
\end{equation}
Using the fact that on $\overline{\{y\in M\mid\chi(y)\neq1\}}$ the function $|df|^2$ is strictly positive, we conclude from
Lemma~\ref{L:F1} that there exists a constant $\epsilon>0$ so that for sufficiently large $u$ and all $w\in\Omega(M;E_u)$ 
satisfying $\supp w\subseteq\overline{\{y\in M\mid\chi(y)\neq1\}}$ we have
\begin{equation}\label{E:B13}
\Re\llangle\Delta_u-\Delta_0)w,w\rrangle\geq\epsilon u^2\|w\|^2.
\end{equation}
Combining \eqref{E:B12} and \eqref{E:B13} the lemma follows easily.
\end{proof}

\begin{lemma}\label{L:3}
There exist constants $u_3'''\geq1$ and $C_3'''\geq0$ so that for any $u\geq u_3'''$ and $w\in\Omega(M;E_u)$
we have
\begin{equation}\label{E:98}
\Re\llangle\Delta_u(\chi w),(1-\chi)w\rrangle \geq -C_3'''\|w\|^2.
\end{equation}
\end{lemma}

\begin{proof}
It suffices to show this inequality for $w'\in\Omega(M;E_u)$ with $\supp w'\subseteq B$.
Indeed, choose a smooth cut-off function $\eta\co M\to[0,1]$ with $\supp\eta\subseteq B$
and $\supp\chi\subseteq\{y\in M\mid \eta(y)=1\}$. Consider $w':=\eta w$. Since $\eta\chi=\chi$
we have $\Delta_u(\chi w')=\Delta_u(\chi w)$, and since $\eta=1$ on $\supp\Delta_u(\chi w)$
we obtain
$$
\llangle\Delta_u(\chi w),(1-\chi)w\rrangle=\llangle\Delta_u(\chi w'),(1-\chi)w'\rrangle.
$$
Moreover, note that $\|w'\|\leq\|w\|$. Hence, if the desired inequality holds for all
$w'\in\Omega(M;E_u)$ with $\supp w'\subseteq B$, it will remain true for arbitrary 
$w\in\Omega(M;E_u)$ with the same constant $C_3'''\geq0$.
In view of Remark~\ref{R:adjoint} the Laplacian $\Delta_u$ coincides with the standard
selfadjoint Witten Laplacian over $B$. For the latter operator this estimate can be found in 
\cite{BFKM96}.
\end{proof}

The above lemmas imply

\begin{proposition}\label{P:CCC}
There exist constants $u_3\geq1$ and $\epsilon_3>0$ so that for all
$u\geq u_3$ and $v'\in V_u^\perp$ we have
$$
\Re\llangle\Delta_uv',v'\rrangle\geq\epsilon_3u\|v'\|^2.
$$
\end{proposition}

\begin{proof}
Suppose $v'\in V_u^\perp$ and write 
$v'=v'_1+v'_2$ with $v_1':=\chi v'$ and $v'_2:=(1-\chi)v'$.
In view of Remark~\ref{R:adjoint} and since $\supp v_1'\subseteq B$ we have
$$
\llangle\Delta_uv_2',v_1'\rrangle
=\llangle v_2',\Delta_uv_1'\rrangle
=\overline{\llangle\Delta_uv_1',v_2'\rrangle}.
$$
Therefore $\Re\llangle\Delta_uv_2',v_1'\rrangle=\Re\llangle\Delta_uv_1',v_2'\rrangle$
and thus
\begin{equation}\label{E:K1}
\Re\llangle\Delta_uv',v'\rrangle=
\Re\llangle\Delta_uv'_1,v'_1\rrangle
+2\Re\llangle\Delta_uv'_1,v'_2\rrangle 
+\Re\llangle\Delta_uv'_2,v'_2\rrangle.
\end{equation}
By Lemma~\ref{L:2} we have
\begin{equation}\label{E:K2}
\Re\llangle\Delta_uv'_1,v'_1\rrangle\geq\epsilon_3'u\|v'_1\|^2,\qquad u\geq u_3'.
\end{equation}
Since $\supp v_2'\subseteq\overline{\{y\in M\mid\chi\neq1\}}$, Lemma~\ref{L:4} implies
\begin{equation}\label{E:K3}
\Re\llangle\Delta_uv_2',v_2'\rrangle\geq\epsilon_3''u^2\|v_2'\|^2,\qquad u\geq u_3''.
\end{equation}
From Lemma~\ref{L:3} we get
\begin{equation}\label{E:K4}
\Re\llangle\Delta_uv_1',v_2'\rrangle\geq-C_3'''\|v'\|^2,\qquad u\geq u_3'''.
\end{equation}
Combining \eqref{E:K1}, \eqref{E:K2}, \eqref{E:K3} and \eqref{E:K4} we obtain 
\begin{equation}\label{E:K5}
\Re\llangle\Delta_uv',v'\rrangle\geq\epsilon_3'u\|v'_1\|^2+\epsilon_3''u^2\|v'_2\|^2-C_3'''\|v'\|^2
\end{equation}
for all $u\geq\max\{u_3',u_3'',u_3'''\}$. Note that adding 
$\|v'\|^2=\|v'_1\|^2+\|v_2'\|^2+2\Re\llangle v_1',v_2'\rrangle$
and $0\leq\|v_1'-v_2'\|^2=\|v_1'\|^2+\|v_2'\|^2-2\Re\llangle v_1',v_2'\rrangle$ yields
$$
\|v'\|^2\leq2\bigl(\|v_1'\|^2+\|v_2'\|^2\bigr).
$$
Combining this with \eqref{E:K5} the statement of the proposition follows immediately.
\end{proof}

\subsection*{The spectral gap}

The estimates in the preceding section permit to establish a widening
gap in the spectrum of $\Delta_u$, as $u\to\infty$. For the precise statement see 
Proposition~\ref{P:specgap} below. This result is a generalization of a well known
``separation of the spectrum property'' in the selfadjoint situation, 
see for instance \cite[Proposition~5.2]{BFKM96}.
We start with the following resolvent estimate.

\begin{lemma}\label{L:DDD} 
There exist constants $u_4\geq1$ and $\epsilon_4>0$ so that for all $u\geq u_4$,
$\lambda\in\C$ and $w\in\Omega(M;E_u)$ we have
$$
\min\bigl\{|\lambda|-e^{-\epsilon_4u},\epsilon_4u-\Re\lambda\bigr\}\|w\|
\leq\|(\Delta_u-\lambda)w\|.
$$
\end{lemma}

\begin{proof}
For $u\geq0$ let us write
$$ 
\pi_u\co\Omega(M;E_u)\to V_u\qquad\text{and}\qquad
\pi_u^\perp\co\Omega(M;E_u)\to V_u^\perp
$$
for the orthogonal projections onto $V_u$ and $V_u^\perp$, respectively.
Let 
$$
\Delta_u=\left(\begin{smallmatrix}\Delta_{1,u}&\Delta_{2,u}\\\Delta_{3,u}&\Delta_{4,u}\end{smallmatrix}\right)
$$
denote the decomposition of $\Delta_u$ with respect to $\Omega(M;E_u)=V_u\oplus V_u^\perp$.
More precisely, we have:
\begin{align*}
\Delta_{1,u}&\co V_u\to V_u&\Delta_{1,u}&:=\pi_u\Delta_u|_{V_u}
\\
\Delta_{2,u}&\co V_u^\perp\to V_u&\Delta_{2,u}&:=\pi_u\Delta_u|_{V_u^\perp}
\\
\Delta_{3,u}&\co V_u\to V_u^\perp&\Delta_{3,u}&:=\pi_u^\perp\Delta_u|_{V_u}
\\
\Delta_{4,u}&\co V_u^\perp\to V_u^\perp&\Delta_{4,u}&:=\pi_u^\perp\Delta_u|_{V_u^\perp}
\end{align*}
Define operators
$$ 
A_u:=\left(\begin{smallmatrix}0&0&\\0&\Delta_{4,u}\end{smallmatrix}\right)
\qquad\text{and}\qquad
B_u:=\left(\begin{smallmatrix}\Delta_{1,u}&\Delta_{2,u}\\\Delta_{3,u}&0\end{smallmatrix}\right).
$$
Clearly, $\Delta_u=A_u+B_u$.
>From Lemma~\ref{L:BBB} we can  see that for $u\geq u_2$ and $w\in\Omega(M;E_u)$ we have
\begin{equation}\label{E:292a}
\|B_uw\|\leq 2e^{-\epsilon_2u}\|w\|.
\end{equation}
Indeed, using $\pi_uB_uw=\pi_u\Delta_uw$ and $\pi_u^\perp B_uw=\pi_u^\perp\Delta_u(\pi_uw)$ we obtain:
\begin{align*}
\|B_uw\|^2
&=\llangle\pi_uB_uw,\pi_uB_uw\rrangle+\llangle\pi_u^\perp B_uw,\pi_u^\perp B_uw\rrangle
\\&=\llangle\Delta_uw,\pi_uB_uw\rrangle+\llangle\Delta_u(\pi_uw),\pi_u^\perp B_uw\rrangle
\\&\leq e^{-\epsilon_2u}\|w\|\|\pi_uB_uw\|+e^{-\epsilon_2u}\|\pi_uw\|\|\pi_u^\perp B_uw\|
\\&\leq2e^{-\epsilon_2u}\|w\|\|B_uw\|
\end{align*}
From Proposition~\ref{P:CCC} and the Cauchy--Schwarz inequality we get
\begin{multline*}
\|(\Delta_{4,u}-\lambda)v'\|\|v'\|
\geq|\llangle(\Delta_{4,u}-\lambda)v',v'\rrangle|
=|\llangle(\Delta_u-\lambda)v',v'\rrangle|
\\\geq\Re\llangle(\Delta_u-\lambda)v',v'\rrangle
=\Re\llangle\Delta_uv',v'\rrangle-\Re\lambda\|v'\|^2
\geq(\epsilon_3u-\Re\lambda)\|v'\|^2
\end{multline*}
and thus
$$
\|(\Delta_{4,u}-\lambda)v'\|\geq(\epsilon_3u-\Re\lambda)\|v'\|
$$
for $u\geq u_3$, all $\lambda\in\C$ and $v'\in V_u^\perp$.
We conclude that
\begin{equation}\label{E:292b}
\|(A_u-\lambda)w\|\geq\min\{|\lambda|,\epsilon_3u-\Re\lambda\}\|w\|
\end{equation}
for $u\geq u_3$, all $\lambda\in\C$ and $w\in\Omega(M;E_u)$.
Combining \eqref{E:292a} and \eqref{E:292b} we find
$$
\|(\Delta_u-\lambda)w\|
\geq\|(A_u-\lambda)w\|-\|B_uw\|
\geq\bigl(\min\{|\lambda|,\epsilon_3u-\Re\lambda\}-2e^{-\epsilon_2u}\bigr)\|w\|
$$
for sufficiently large $u$, all $\lambda\in\C$ and $w\in\Omega(M;E_u)$. 
The statement now follows with an appropriate choice of $\epsilon_4$ and $u_4$.
\end{proof}

\begin{proposition}\label{P:specgap}
Let $\epsilon_4>0$ be the constant from Lemma~\ref{L:DDD}. Then, 
for sufficiently large $u$, we have
$$
\Spec(\Delta_u)\subseteq\bigl\{\lambda\in\C\bigm| |\lambda|\leq e^{-\epsilon_4u}\bigr\}
\cup\bigl\{\lambda\in\C\bigm| \Re\lambda\geq\epsilon_4u\bigr\}.
$$
\end{proposition}

\begin{proof}
Suppose $\lambda\in\Spec(\Delta_u)$. Then there exists $w\in\Omega(M;E_u)$ with
$(\Delta_u-\lambda)w=0$ and $\|w\|\neq0$. Assuming $u$ is sufficiently large, see 
Lemma~\ref{L:DDD}, we conclude
$\min\bigl\{|\lambda|-e^{-\epsilon_4u},\epsilon_4u-\Re\lambda\bigr\}\leq0$, and the statement follows.
\end{proof}

\subsection*{Approximation of the small complex}

We will now show that, as $u\to\infty$, the subspace $V_u\subseteq\Omega(M;E_u)$
approximates $\Omega_\sm(M;E_u)$ well with respect to every $C^k$--norm
on $\Omega(M;E_u)$. For the precise statements 
see Proposition~\ref{P:QQQ} and Proposition~\ref{P:ontoness} below.
Recall that $\Omega_\sm(M;E_u)$ denotes the sum of eigen spaces of $\Delta_u$ whose corresponding eigen values
have real part at most $1$.

For $s\in\N$ and $w\in\Omega(M;E_u)$ let $\|w\|_s$ denote the (a fixed) 
Sobolev $s$-norm. We continue to write $\|w\|=\|w\|_0$.
We will start with the following improvement of the resolvent estimate in
Lemma~\ref{L:DDD}.

\begin{lemma}\label{L:TTT}
Let $u_4\geq1$ and $\epsilon_4>0$ be the constants from Lemma~\ref{L:DDD}.
For every $s\in\N$ there exists a constant $C_{4,s}\geq0$ so that for
all $u\geq u_4$, $\lambda\in\C$ and $w\in\Omega(M;E_u)$ we have
\begin{equation}\label{E:314}
\min\bigl\{|\lambda|-e^{-\epsilon_4u},\epsilon_4u-\Re\lambda\bigr\}\|w\|_s
\leq C_{4,s}(u^2+|\lambda|)^s\|(\Delta_u-\lambda)w\|_s.
\end{equation}
Moreover, there exists a constant $\tilde C_{4,2s}\geq0$ so that for all 
$u\geq u_4$, $\lambda\in\C$ and $\omega\in\Omega(M;E_u)$ we even have
\begin{equation}\label{E:314b}
\min\bigl\{|\lambda|-e^{-\epsilon_4u},\epsilon_4u-\Re\lambda\bigr\}\|w\|_{2s}
\leq\tilde C_{4,2s}(u^2+|\lambda|)^s\|(\Delta_u-\lambda)w\|_{2s}.
\end{equation}
\end{lemma}

\begin{proof}
We will construct the constants $C_{4,s}\geq0$ by induction on $s$.
For $s=0$ the statement was proved in Lemma~\ref{L:DDD}.
The induction is based on an argument that is used in the selfadjoint situation too, see
\cite[proof of Proposition~5.4]{BFKM96} or \cite[proof of Theorem~8.8]{BZ92}.

From the ellipticity of $\Delta_0$ we get constants $C_s'\geq0$
so that for all $w\in\Omega(M;E_u)$ 
\begin{equation}\label{E:A01}
\|w\|_{s+1}\leq\|w\|_{s+2}\leq C'_s\bigl(\|\Delta_0w\|_s+\|w\|_s\bigr).
\end{equation}
From Lemma~\ref{L:F1} we obtain constants $C_s''\geq0$ such that 
for all $u\geq u_4$ and $w\in\Omega(M;E_u)$ 
\begin{equation}\label{E:A02}
\|(\Delta_u-\Delta_0)w\|_s\leq C_s''u^2\|w\|_s.
\end{equation}
Combining \eqref{E:A01} and \eqref{E:A02} we obtain constants $C_s'''\geq0$
so that for all $u\geq u_4$, $\lambda\in\C$ and $w\in\Omega(M;E_u)$ we have
\begin{equation}\label{E:A03}
\|w\|_{s+1}\leq C'''_s
\Bigl(\|(\Delta_u-\lambda)w\|_s+(u^2+|\lambda|)\|w\|_s\Bigr).
\end{equation}
By induction we may assume that there exists a constant $C_{4,s}\geq0$ so that 
\eqref{E:314} holds. Combining \eqref{E:314} with \eqref{E:A03} 
and using $\min\bigl\{|\lambda|-e^{-\epsilon_4u},\epsilon_4u-\Re\lambda\bigr\}\leq|\lambda|$
we find a constant $C_{4,s+1}\geq0$ so that for all $u\geq u_4$, $\lambda\in\C$ and
$\omega\in\Omega(M;E_u)$
$$
\min\bigl\{|\lambda|-e^{-\epsilon_4u},\epsilon_4u-\Re\lambda\bigr\}\|w\|_{s+1}
\leq C_{4,s+1}(u^2+|\lambda|)^{s+1}\|(\Delta_u-\lambda)w\|_s
$$ 
Now use $\|(\Delta_u-\lambda)w\|_s\leq\|(\Delta_u-\lambda)w\|_{s+1}$ to complete the induction.
The proof of \eqref{E:314b} is similar.
\end{proof}

Let $Q_u\co\Omega(M;E_u)\to\Omega_\sm(M;E_u)$ denote the spectral projection.

\begin{lemma}\label{L:SSS} 
For every $s\in\N$ there exist constants $u_{5,s}\geq1$ and $\epsilon_{5,s}>0$ so that
for all $u\geq u_{5,s}$ and $v\in V_u$ we have 
$$
\|Q_uv-v\|_s\leq e^{-\epsilon_{5,s}u}\|v\|.
$$
\end{lemma}

\begin{proof}
By Proposition~\ref{P:specgap}, for sufficiently large 
$u$, 
we have $\Spec(\Delta_u) \cap S^1=\emptyset$
and hence $Q_u$ is given by the Riesz projector
$$
Q_u=\frac1{2\pi\ii}\int_{S^1}(\lambda-\Delta_u)^{-1}d\lambda.
$$
Since 
$(\lambda-\Delta_u)^{-1}-\lambda^{-1}=\lambda^{-1}(\lambda-\Delta_u)^{-1}\Delta_u$
we obtain, for $v\in\Omega(M;E_u)$,
\begin{equation}\label{E:A04}
Q_uv-v=\frac1{2\pi\ii}\int_{S^1}\lambda^{-1}(\lambda-\Delta_u)^{-1}\Delta_uv\,d\lambda.
\end{equation}
From Lemma~\ref{L:TTT} and Lemma~\ref{L:AAA} we easily infer the existence of constants 
$u_{5,s}\geq1$ and $\epsilon_{5,s}>0$ so that for all $s\in\N$, 
$u\geq u_{5,s}$, $\lambda\in S^1$ and $v\in V_u$ we have
$$
\|\lambda^{-1}(\lambda-\Delta_u)^{-1}\Delta_uv\|_s\leq e^{-\epsilon_{5,s}u}\|v\|.
$$
The lemma now follows easily by combining this estimate with \eqref{E:A04}.
\end{proof}

\begin{proposition}\label{P:QQQ}
There exist constants $u_{6,k}\geq1$ and $\epsilon_{6,k}>0$ so that for 
$k\in\N$, $u\geq u_{6,k}$ and $v\in V_u$ we have
$$
\|Q_uv-v\|_{C^k}\leq e^{-\epsilon_{6,k}u}\|v\|.
$$
\end{proposition}

\begin{proof}
This follows from Lemma~\ref{L:SSS} and the Sobolev embedding theorem.
\end{proof}

In the selfadjoint case the following estimate is an immediate consequence of Proposition~\ref{P:specgap}.
In our situation we will have to use the resolvent estimate from Lemma~\ref{L:TTT}.

\begin{lemma}\label{L:RRR}
For every $s\in\N$ there exist constants $u_{7,s}\geq1$ and $\epsilon_{7,s}>0$ so that for all
$u\geq u_{7,s}$ and all $w\in\Omega_\sm(M;E_u)$ we have
$$
\|\Delta_uw\|_s\leq e^{-\epsilon_{7,s}u}\|w\|_s.
$$
\end{lemma}

\begin{proof}
Let $\epsilon_4>0$ be the constant from Lemma~\ref{L:DDD}, and set $\rho_u:=2e^{-\epsilon_4u}$.
Assume $u\geq0$ is sufficiently large so that all eigen values with real part at most $1$ are contained in the interior of the 
circle $\rho_uS^1$ of radius $\rho_u$, see Proposition~\ref{P:specgap}. Then
\begin{equation}\label{E:aux1}
\Delta_uQ_u=Q_u\Delta_u=\frac1{2\pi\ii}\int_{\rho_uS^1}\lambda(\lambda-\Delta_u)^{-1}d\lambda.
\end{equation}
For $\lambda\in\rho_uS^1$ and sufficiently large $u$, Lemma~\ref{L:TTT}, see \eqref{E:314}, provides the estimate
$$
\|\lambda(\lambda-\Delta_u)^{-1}w\|_s\leq2C_{4,s}(u^2+2)^s\|w\|_s,\qquad w\in\Omega(M;E_u).
$$
Combining this with \eqref{E:aux1} we obtain, for $w\in\Omega(M;E_u)$,
$$
\|\Delta_uQ_uw\|_s\leq2C_{4,s}(u^2+2)^s\rho_u\|w\|_s.
$$
Choosing $\epsilon_{7,s}>0$ and $u_{7,s}\geq1$ appropriately we get
$$
\|\Delta_uQ_uw\|_s\leq e^{-\epsilon_{7,s}u}\|w\|_s
$$
for all $u\geq u_{7,s}$ and $w\in\Omega(M;E_u)$. The statement follows immediately.
\end{proof}

\begin{proposition}\label{P:ontoness}
For sufficiently large $u$ the restriction of the spectral projection
$Q_u\co V_u\to\Omega_\sm(M;E_u)$ is an isomorphism.
\end{proposition}

\begin{proof} 
Recall that $\Omega_\sm(M;E_u)$ is the image of the spectral projection $Q_u$, 
and that $\beta(Q_uv,w)=\beta(v,Q_uw)$ for all $v,w\in\Omega(M;E_u)$. Hence,
\begin{equation}\label{E:655}
\Omega_\sm(M;E_u)\cap(Q_uV_u)^{\perp_\beta}=\Omega_\sm(M;E_u)\cap V_u^\perp.
\end{equation}
Here $(Q_uV_u)^{\perp_\beta}$ denotes the $\beta$-orthogonal complement of $Q_uV_u$,
but see also Observation~\ref{O:Vubar}.
It follows from Proposition~\ref{P:CCC} and Lemma~\ref{L:RRR} 
that $\Omega_\sm(M;E_u)\cap V_u^\perp=0$ for sufficiently large $u$.
Together with \eqref{E:655} this shows that $Q_uV_u=\Omega_\sm(M;E_u)$, and
therefore $Q_u\co V_u\to\Omega_\sm(M;E_u)$ is onto, for sufficiently large $u$.
The injectivity follows from Lemma~\ref{L:SSS} with $s=0$.
\end{proof}

\subsection*{The integration on the small complex}

For the sake of notational simplicity let us introduce the notation, $u>0$,
\begin{align*}
I_u:=\eta_u{\Int_u}|_{V_u}&\co V_u\to C(X;E_u)
\\
I_{u,\sm}:=\eta_u\Int_{u,\sm}&\co\Omega_\sm(M;E_u)\to C(X;E_u)
\end{align*}
where $\eta_u\co C(X;E_u)\to C(X;E_u)$ denotes the scaling isomorphism
introduced at the beginning of this section, see \eqref{E:scaling}.
Moreover, we will write $\beta_u:=\beta|_{V_u}$ for the restriction of the 
bilinear form $\beta$ to $V_u$, and we will continue to write
$\beta_{u,\sm}=\beta|_{\Omega_\sm(M;E_u)}$ for the restriction of the bilinear 
form $\beta$ to $\Omega_\sm(M;E_u)$. Recall that we write
$\langle\cdot,\cdot\rangle_{\mathcal X}$ and $|\cdot|_{\mathcal X}$ for the
Hermitian inner product and its associated norm on $C(X;E_u)$.
Finally, we recall that $b_{\mathcal X}$ denotes the bilinear form on
$C(X;E_u)$, see \eqref{E:A50}.

\begin{lemma}\label{L:I1}
For sufficiently large $u$ the mapping $I_u\co V_u\to C(X;E_u)$
is an isomorphism. Moreover, there exists a constant $\epsilon_8>0$ so that, 
as $u\to\infty$,
\begin{equation}\label{E:A99}
(I_u)_*\llangle\cdot,\cdot\rrangle=\langle\cdot,\cdot\rangle_{\mathcal X}
+\LO(e^{-\epsilon_8u})
\end{equation}
and
\begin{equation}\label{E:A98}
(I_u)_*\beta_u=b_{\mathcal X}+\LO(e^{-\epsilon_8u})
\end{equation}
\end{lemma}

\begin{proof}
We claim that there exist constants $u_8\geq1$ and $\epsilon_8>0$ so that
for all $u\geq u_8$ and $v_1,v_2\in V_u$ we have
\begin{equation}\label{E:A20}
\Bigl|\langle I_uv_1,I_uv_2\rangle_{\mathcal X}-\llangle v_1,v_2\rrangle\Bigr|
\leq e^{-\epsilon_8u}|I_uv_1|_{\mathcal X}|I_uv_2|_{\mathcal X}.
\end{equation}
If such  estimate is established we immediately see that $I_u\co V_u\to C(X;E_u)$
is injective for $u\geq u_8$. Since $V_u$ and $C(X;E_u)$
obviously have the same dimension, the first assertion of the lemma follows.
>From \eqref{E:A20} we obtain, for $u\geq u_8$ and $a_1,a_2\in C(X;E_u)$,
$$
\Bigl|\langle a_1,a_2\rangle_{\mathcal X}-\llangle I_u^{-1}a_1,I_u^{-1}a_2\rrangle\Bigr|
\leq e^{-\epsilon_8u}|a_1|_{\mathcal X}|a_2|_{\mathcal X}.
$$ 
from which we infer \eqref{E:A99}. Combining the last equation with
\eqref{E:A50}, \eqref{E:A51}, \eqref{EEEE:2} and using the fact that 
$I_u\co V_u\to C(X;E_u)$ intertwines the complex conjugations we also obtain,
for $a_1,a_2\in C(X;E_u)$,
$$ 
\Bigl|b_{\mathcal X}(a_1,a_2)-\beta_u(I_u^{-1}a_1,I_u^{-1}a_2)\Bigr|
\leq e^{-\epsilon_8u}|a_1|_{\mathcal X}|a_2|_{\mathcal X}
$$
and thus \eqref{E:A98}.

It thus remains to establish the estimate \eqref{E:A20}. 
To do so, suppose $x\in\mathcal X$ and set $q:=\ind(x)$. For $e_1,e_2\in E_x$ we have, see \eqref{E:phiuxe},
$$
\llangle\phi_{u,e_1},\phi_{u,e_2}\rrangle
=\int_{B_x}(\sigma\circ r_x)^2e^{-ur_x^2}\langle e_1,e_2\rangle\vvol_g
=\int_{\R^n}\sigma(|z|)^2e^{-u|z|^2}dz\langle e_1,e_2\rangle
$$
and thus there exist constants $C'\geq0$ and $\epsilon'>0$ so that for all $u\geq1$
and $e_1,e_2\in E_x$
\begin{equation}\label{E:B20}
\Bigl|\llangle\phi_{u,e_1},\phi_{u,e_2}\rrangle-(\pi/u)^{n/2}\langle e_1,e_2\rangle\Bigr|\leq
C'e^{-u\epsilon'}|e_1||e_2|
\end{equation}
Further, keeping \eqref{E:SmB}, \eqref{E:XX2} and \eqref{E:M1} in mind, we have for $e_1,e_2\in E_x$,
\begin{align*} 
\bigl\langle I_u\phi_{u,e_1},&I_u\phi_{u,e_2}\bigr\rangle_{\mathcal X}
\\
&=\biggl((\pi/u)^{n/4-q/2}\int_{W_x^-\cap B_x}\!\!\!\!\!\!\!
(\sigma\circ r_x)e^{-ur_x^2/2}e^{u(f-f(x))}\Omega_x^-\biggr)^2\langle e_1,e_2\rangle
\\
&=(\pi/u)^{n/2-q}\biggl(\int_{\R^q}\sigma(|z|)
e^{-u|z|^2}dz\biggr)^2\langle e_1,e_2\rangle
\end{align*}
and hence there exist constants $C''\geq0$ and $\epsilon''>0$ so that for all $u\geq1$
and $e_1,e_2\in E_x$
\begin{equation}\label{E:B21}
\Bigl|\bigl\langle I_u\phi_{u,e_1},I_u\phi_{u,e_2}\bigr\rangle_{\mathcal X}
-(\pi/u)^{n/2}\langle e_1,e_2\rangle\Bigr|
\leq C''e^{-u\epsilon''}|e_1||e_2|.
\end{equation}
Note that this estimate also implies, for $e\in E_x$ and $u\geq1$
\begin{equation}\label{E:B22}
\Bigl((\pi/u)^{n/2}-C''e^{-u\epsilon''}\Bigr)|e|^2\leq|I_u\phi_{u,e}|_{\mathcal X}^2.
\end{equation}

Combining \eqref{E:B20}, \eqref{E:B21} and \eqref{E:B22} we find constants
$u_8\geq1$ and $\epsilon_8>0$ so that for all $u\geq u_8$, $x\in\mathcal X$,
and $e_1,e_2\in E_x$ we have
$$
\Bigl|
\bigl\langle I_u\phi_{u,e_1},I_u\phi_{u,e_2}\bigr\rangle_{\mathcal X}
-\llangle\phi_{u,e_1},\phi_{u,e_2}\rrangle
\Bigr| 
\leq e^{-u\epsilon_8}|I_u\phi_{u,e_1}|_{\mathcal X}|I_u\phi_{u,e_2}|_{\mathcal X}.
$$
For $x,y\in\mathcal X$, $x\neq y$, $e_1\in E_x$, $e_2\in E_y$ we clearly have 
$\langle I_u\phi_{u,e_1},I_u\phi_{u,e_2}\rangle=0=\llangle\phi_{u,e_1},\phi_{u,e_2}\rrangle$, see \eqref{E:SmB}.
The estimate \eqref{E:A20} now follows easily, see Definition~\ref{D:Vu}, and the proof is complete.
\end{proof}

\begin{lemma}\label{L:YYY}
There exist constants $u_9\geq1$ and $\epsilon_9>0$ so that for all $u\geq u_9$
and $v\in V_u$ we have
$$
\bigl|I_{u,\sm}Q_uv-I_uv\bigr|_{\mathcal X}\leq e^{-\epsilon_9u}\|v\|.
$$
\end{lemma}

\begin{proof}
There exists a constant $C\geq0$ so that for all
$u\geq0$ and $w\in\Omega(M;E_u)$ 
$$
|\Int_uw|_{\mathcal X}\leq C\|w\|_{C^0}.
$$
To see this we use the compactification result for the unstable manifolds. The
uniformity in $u$ is guaranteed by the 
relation $\Int_uw=e^{-uf}\Int_0(e^{uf}w)$ and the fact that the function
$e^{uf}$ restricted to an unstable manifold $W_x^-$ will attain its maximum
at the critical point $x\in\mathcal X$. The statement then follows from
Proposition~\ref{P:QQQ} with $k=0$.
\end{proof}

\begin{lemma}\label{L:XXX}
There exist constant $u_{10}\geq1$ and $\epsilon_{10}>0$ so that
for all $u\geq u_{10}$ and $v_1,v_2\in V_u$ we have
$$
\bigl|\beta_{u,\sm}(Q_uv_1,Q_uv_2)-\beta_u(v_1,v_2)\bigr|\leq e^{-\epsilon_{10}u}\|v_1\|\|v_2\|.
$$
\end{lemma}

\begin{proof}
From \eqref{EEEE:2}, \eqref{EEEE:4} and the Cauchy--Schwarz inequality we obtain
$$
|\beta(w_1,w_2)|\leq\|w_1\|\|w_2\|
$$ 
for all $w_1,w_2\in\Omega(M;E_u)$. Hence:
\begin{align*}
\bigl|\beta(Q_uw_1,&Q_uw_2)-\beta(w_1,w_2)\bigr|
\\&\leq\bigl|\beta(Q_uw_1-w_1,Q_uw_2)\bigr|+\bigl|\beta(w_1,Q_uw_2-w_2)\bigr|
\\&\leq\|Q_uw_1-w_1\|\|Q_uw_2\|+\|w_1\|\|Q_uw_2-w_2\|
\\&\leq\|Q_uw_1-w_1\|\bigl(\|Q_uw_2-w_2\|+\|w_2\|\bigr)+\|w_1\|\|Q_uw_2-w_2\|
\end{align*}
The statement thus follows from Lemma~\ref{L:SSS} with $s=0$.
\end{proof}

We are now in the position to put the pieces together and provide a

\begin{proof}[Proof of Theorem~\ref{T:WHS}]
From \eqref{E:A99} we obtain a constant $C'\geq0$ so that for sufficiently
large $u$ and $a\in C(X;E_u)$ we have
\begin{equation}\label{E:A30}
\|I_u^{-1}a\|\leq C'|a|_{\mathcal X}.
\end{equation}
Combining this with Lemma~\ref{L:YYY} we obtain constant $\epsilon'>0$ so
that, as $u\to\infty$,
\begin{equation}\label{E:A97}
I_{u,\sm}Q_uI_u^{-1}=\id_{C(X;E_u)}+\LO(e^{-\epsilon'u}).
\end{equation}
Particularly, the mapping $I_{u,\sm}Q_uI_u^{-1}\co C(X;E_u)\to C(X;E_u)$
is an isomorphism, for sufficiently large $u$. 
Hence $I_{u,\sm}\co\Omega_\sm(M;E_u)\to C(X;E_u)$ is an isomorphism 
for sufficiently large $u$, see Proposition~\ref{P:ontoness} and Lemma~\ref{L:I1}.
Since the scaling $\eta_u$ is an isomorphism for
all $u>0$, we see that $\Int_{u,\sm}=\eta_u^{-1}I_{u,\sm}\co\Omega_\sm(M;E_u)\to
C(X;E_u)$ is an isomorphism too.
This proves the first claim of Theorem~\ref{T:WHS}.

Next note that \eqref{E:A30} and \eqref{E:A97} provide a constant
$C''\geq0$ so that for sufficiently large $u$ and $a\in C(X;E_u)$ we have
\begin{equation}\label{E:A31}
\|(I_{u,\sm}Q_u)^{-1}a\|\leq C''|a|_{\mathcal X}.
\end{equation}
Combining this with Lemma~\ref{L:XXX} we find a constant $\epsilon''>0$ so that
\begin{equation}\label{E:A41}
(I_{u,\sm}Q_u)_*\bigl(Q_u^*\beta_{u,\sm}-\beta_u\bigr)
=\LO(e^{-\epsilon'' u}).
\end{equation}
From \eqref{E:A97} and \eqref{E:A98} we find a constant $\epsilon'''>0$ so that, as
$u\to\infty$,
\begin{equation}\label{E:A40}
(I_{u,\sm}Q_u)_*\beta_u
=(I_{u,\sm}Q_uI_u^{-1})_*(I_u)_*\beta_u
=b_{\mathcal X}+\LO(e^{-\epsilon'''u}).
\end{equation}
Clearly, \eqref{E:A41} and \eqref{E:A40} imply the existence 
of a constant $\epsilon>0$ so that, as $u\to\infty$,
$$
(I_{u,\sm})_*\beta_{u,\sm}
=(I_{u,\sm}Q_u)_*(Q_u^*\beta_{u,\sm}-\beta_u)+(I_{u,\sm}Q_u)_*\beta_u
=b_{\mathcal X}+\LO(e^{-\epsilon u}).
$$
This completes the proof of Theorem~\ref{T:WHS}.
\end{proof}

\subsection*{A uniform estimate for the heat trace}

The aim of this section is to establish Theorem~\ref{T:sesti} below. This estimate 
generalizes 
\cite[Theorem~7.7~and Theorem ~7.8]{BZ92} to the non-selfadjoint situation. Similar estimates can be found in \cite{BFK96}.
We will proceed in the spirit of \cite{BZ92}, see also \cite[Section~9]{BL91}.

Let $P_u:=1-Q_u$ denote the spectral projection onto the sum of eigen spaces whose corresponding eigen values
have real part larger than $1$. Moreover, for a trace class operator $A$, let $\|A\|_{\tr}:=\tr(A^*A)$ denote 
the trace norm, see \cite[Appendix~A.3.4]{S01}.

\begin{theorem}\label{T:sesti}
There exist constants $\epsilon>0$ and $u_0\geq1$ with the following property. For every $p>n/2$, $p\in\N$, 
there exists a constant $C_p\geq0$ so that for all $t>0$ and $u\geq u_0$ we have
\begin{equation}\label{E:A16}
\|\exp(-t\Delta_u)P_u\|_{\tr}\leq C_pe^{-\epsilon tu}\Bigl(\frac{u^{p(p+1)+2}}{t^{p-1}}+\frac{u}{t^p}\Bigr).
\end{equation}
\end{theorem}
The proof of Theorem~\ref{T:sesti} is based on estimates for the resolvent of $\Delta_u$ which are uniform in $u$, see
Proposition~\ref{P:sestia} and \ref{P:sestib} below. Proposition~\ref{P:sestia} is of very general nature, whereas
Proposition~\ref{P:sestib} makes use of the Witten--Helffer--Sj\"ostrand estimates in an essential way.
We fix an angle $0<\theta<\pi/4$. For an operator $A$ and $s_1,s_2\in\N$ we will write 
$\|A\|_{s_1,s_2}=\sup_{w\neq0}\|Aw\|_{s_2}/\|w\|_{s_1}$, where $\|w\|_s$  denotes the Sobolev $s$-norm of $w$.

\begin{proposition}\label{P:sestia}
For every $p\in\N$ there exists constants $\alpha_p\geq0$ and
$C_p'\geq0$ so that the following holds.
If $u\geq 1$, $\lambda\in\C$, $|\arg\lambda|\geq\theta$ and $\alpha_pu^2\leq|\lambda|$ then
$\lambda\notin\Spec(\Delta_u)$ and 
$$
\|(\Delta_u-\lambda)^{-p}\|_{0,2p}\leq C_p'.
$$
\end{proposition}

For the proof of Proposition~\ref{P:sestia} we need the following two lemmas.

\begin{lemma}\label{L:sestia}
For every $s\in\N$ there exists a constant $\tilde C_s\geq0$ with the following property. If $u\geq1$ and
$\lambda\notin\Spec(\Delta_u)$ then
$$
\|(\Delta_u-\lambda)^{-1}\|_{s,s+2}\leq\tilde C_s\Bigl(1+(|\lambda|+u^2)\|(\Delta_u-\lambda)^{-1}\|_{s,s}\Bigr).
$$
\end{lemma}

\begin{proof}
By ellipticity there exists a constant $\bar C_s\geq0$ so that for every $w\in\Omega(M;E_u)$ we have:
\begin{align*}
\|w\|_{s+2}&\leq\bar C_s\bigl(\|\Delta_0w\|_s+\|w\|_s\bigr)
\\&=\bar C_s\Bigl(\|(\Delta_u-\lambda+\lambda-\Delta_u+\Delta_0)(w)\|_s+\|w\|_s\Bigr)
\\&\leq\bar C_s\Bigl(\|(\Delta_u-\lambda)w\|_s+\bigl(1+|\lambda|+\|\Delta_u-\Delta_0\|_{s,s}\bigr)\|w\|_s\Bigr)
\\&\leq\bar C_s\Bigl(1+\bigl(1+|\lambda|+\|\Delta_u-\Delta_0\|_{s,s}\bigr)\|(\Delta_u-\lambda)^{-1}\|_{s,s}\Bigr)\|(\Delta_u-\lambda)w\|_s
\end{align*} 
Since $\|\Delta_u-\Delta_0\|_{s,s}=O(u^2)$ we find a constant $\tilde C_s\geq0$ so that
$$
\|w\|_{s+2}\leq\tilde C_s\Bigl(1+(|\lambda|+u^2)\|(\Delta_u-\lambda)^{-1}\|_{s,s}\Bigr)\|(\Delta_u-\lambda)w\|_s,
$$
and this implies the statement of the lemma.
\end{proof}

\begin{lemma}\label{L:sestib}
For every $s\in\N$ there exist constants $\tilde\alpha_s\geq0$ and $\tilde C_s'\geq0$ with the following property. For all $u\geq1$
and $\lambda\in\C$ with $\tilde\alpha_su^2\leq|\lambda|$ and $|\arg\lambda|\geq\theta$ we have $\lambda\notin\Spec(\Delta_u)$ and
$$
\|(\Delta_u-\lambda)^{-1}\|_{s,s}\leq\frac{\tilde C_s'}{|\lambda|}.
$$
\end{lemma}

\begin{proof}
In view of \cite[Corollary~9.2]{S01} there exists constants $R\geq0$ and $\tilde C_s'>0$ so that the following holds.
If $|\lambda|\geq R$ and $|\arg\lambda|\geq\theta$ then $\lambda\notin\Spec(\Delta_0)$ and
\begin{equation}\label{E:A7}
\|(\Delta_0-\lambda)^{-1}\|_{s,s}\leq\frac{\tilde C_s'}{2|\lambda|}.
\end{equation}
Choose $\tilde\alpha_s\geq R$ such that for all $u\geq1$ we have
\begin{equation}\label{E:A88}
\|\Delta_u-\Delta_0\|_{s,s}\leq\frac{\tilde\alpha_s}{\tilde C_s'}u^2.
\end{equation}
If $u\geq1$, $\tilde\alpha_su^2\leq|\lambda|$, $|\arg\lambda|\geq\theta$ and $\omega\in\Omega(M;E_u)$, then
combining \eqref{E:A7} and \eqref{E:A88} we obtain
\begin{equation*}
\begin{aligned}
\|(\Delta_u-\lambda)w\|_s
\geq\|(\Delta_0-\lambda)w\|_s-\|(\Delta_u-\Delta_0)w\|_s\\
\geq\Bigl(\frac{2|\lambda|}{\tilde C_s'}-\frac{\tilde\alpha_su^2}{\tilde C_s'}\Bigr)\|w\|_s
\geq\frac{|\lambda|}{\tilde C_s'}\|w\|_s
\end{aligned}
\end{equation*}
and the lemma follows.
\end{proof}

\begin{proof}[Proof of Proposition~\ref{P:sestia}]
For $s\in\N$ set $\tilde C_s'':=\tilde C_s(1+(1+\frac1{\tilde\alpha_s})\tilde C_s')$ where $\tilde C_s$, $\tilde\alpha_s$ and $\tilde C_s'$ 
are the constants from Lemma~\ref{L:sestia} and \ref{L:sestib}. Then
$$
\|(\Delta_u-\lambda)^{-1}\|_{s,s+2}\leq\tilde C_s''
$$
for all $u\geq1$, $\tilde\alpha_su^2\leq|\lambda|$  and $|\arg\lambda|\geq\theta$.
Set $\alpha_p:=\max\bigl\{\tilde\alpha_0,\tilde\alpha_2,\tilde,\dotsc,\tilde\alpha_{2p-2}\bigr\}$.
Then
$$
\|(\Delta_u-\lambda)^{-p}\|_{0,2p}
\leq\prod_{s=0}^{p-1}\|(\Delta_u-\lambda)^{-1}\|_{2s,2s+2}
\leq\prod_{s=0}^{p-1}\tilde C_s''
$$
for all $u\geq1$, $\alpha_pu^2\leq|\lambda|$ and $|\arg\lambda|\geq\theta$. The proposition now follows
with $C_p':=\prod_{s=0}^{p-1}\tilde C_s''$.
\end{proof}

Let us introduce the notation:
\begin{align*}
\Spec_\sm(\Delta_u)&:=\bigl\{\lambda\in\Spec(\Delta_u)\bigm|\Re\lambda\leq1\bigr\}
\\
\Spec_\la(\Delta_u)&:=\bigl\{\lambda\in\Spec(\Delta_u)\bigm|\Re\lambda>1\bigr\}
\end{align*}

\begin{proposition}\label{P:sestib}
There exist constants $\epsilon>0$, $u_0\geq1$ and for every $p\in\N$ a constant $C_p''\geq0$ so that the following holds.
If $u\geq u_0$ then  (by Proposition~\ref{P:specgap})
\begin{align}
\label{E:sestiSa}
\Spec_\sm(\Delta_u)&\subseteq\bigl\{\lambda\in\C\bigm|\Re\lambda<\epsilon u\bigr\},
\\\label{E:sestiSb}
\Spec_\la(\Delta_u)&\subseteq\bigl\{\lambda\in\C\bigm|\Re\lambda>2\epsilon u\bigr\},
\end{align}
and for $\lambda\in\C$ with
$\epsilon u\leq\Re\lambda\leq2\epsilon u$, we have
$$ 
\|(\Delta_u-\lambda)^{-p}\|_{0,2p}\leq C_p''(|\lambda|+u^2)^{p(p+1)/2}.
$$
\end{proposition}

\begin{proof}[Proof of Proposition~\ref{P:sestib}]
Set $\epsilon:=\epsilon_4/4$ where $\epsilon_4>0$ denotes the constant in Proposition~\ref{P:specgap}, and choose $u_0\geq1$ sufficiently large
so that for $u\geq u_0$ \eqref{E:sestiSa} and \eqref{E:sestiSb} hold. 
In view of Lemma~\ref{L:TTT}, see \eqref{E:314b}, we can increse $u_0$ so that for $u\geq u_0$, $\epsilon u\leq\Re\lambda\leq2\epsilon u$ we have
$$
\|(\Delta_u-\lambda)^{-1}\|_{2s,2s}\leq(u^2+|\lambda|)^s,\qquad s=0,1,\dotsc,p-1.
$$
Then, using Lemma~\ref{L:sestia}, 
$$
\|(\Delta_u-\lambda)^{-1}\|_{2s,2s+2}
\leq\tilde C_{2s}\bigl(1+(|\lambda|+u^2)^{s+1}\bigr)
\leq 2\tilde C_{2s}(|\lambda|+u^2)^{s+1}
$$
and we obtain
$$
\|(\Delta_u-\lambda)^{-p}\|_{0,2p}
\leq\prod_{s=0}^{p-1}\|(\Delta_u-\lambda)^{-1}\|_{2s,2s+2}
\leq(|\lambda|+u^2)^{p(p+1)/2}\prod_{s=0}^{p-1}2\tilde C_{2s}.
$$
The proposition thus follows with $C_p'':=\prod_{s=0}^{p-1}2\tilde C_{2s}$.
\end{proof}

We will now use Propositions~\ref{P:sestia} and \ref{P:sestib} to provide a

\begin{proof}[Proof of Theorem~\ref{T:sesti}]
We are going to use the constants $\alpha_p\geq0$, $C_p'\geq0$, $\epsilon>0$, $u_0\geq1$ and $C_p''\geq0$
from Propositions~\ref{P:sestia} and \ref{P:sestib}. Increasing $\alpha_p$, we may 
assume that $\alpha_p\geq\epsilon$. For $u\geq u_0$ we consider the contour\footnote{The angular contour bisected 
by the positive real axis.} $\Gamma_u$ parametrized by
$$ 
\lambda_u\co\R\to\C,\qquad
\lambda_u(x):=\epsilon u+\epsilon u|x|-\alpha_pu^2x\ii.
$$
Note that for $u\geq u_0$, $t>0$ and $x\in\R$ we have
\begin{equation}\label{E:A15}
|e^{-t\lambda_u(x)}|=e^{-\epsilon tu}e^{-\epsilon tu|x|}
\qquad\text{and}\qquad
|\lambda_u'(x)|\leq \sqrt2\alpha_pu^2.
\end{equation}
Observe that for $|x|\leq1$ we have $\epsilon u\leq\Re\lambda_u(x)\leq2\epsilon u$ and $|\lambda_u(x)|\leq3\alpha_pu^2$, hence
Proposition~\ref{P:sestib} tells that
\begin{equation}\label{E:A13}
\|(\lambda_u(x)-\Delta_u)^{-p}\|_{0,2p}\leq C''_p(3\alpha_p+1)u^{p(p+1)}
\end{equation}
for all $|x|\leq1$ and $u\geq u_0$.
Moreover, if $|x|\geq1$ then $\alpha_pu^2\leq|\lambda_u(x)|$ and $|\arg\lambda_u(x)|\geq\pi/4\geq\theta$, hence Proposition~\ref{P:sestia} yields
\begin{equation}\label{E:A10}
\|(\lambda_u(x)-\Delta_u)^{-p}\|_{0,2p}\leq C_p'
\end{equation}
for all $|x|\geq1$ and $u\geq u_0$. Further we have:
\begin{align*}
\exp(-t\Delta_u)P_u
&=
\frac1{2\pi\ii}\int_{\Gamma_u}e^{-t\lambda}(\lambda-\Delta_u)^{-1}\ d\lambda
\\&=
\frac{(p-1)!}{2\pi\ii}\frac{(-1)^{p-1}}{t^{p-1}}\int_{\Gamma_u}e^{-t\lambda}(\lambda-\Delta_u)^{-p}\ d\lambda
\\&=
\frac{(p-1)!}{2\pi\ii}\frac{(-1)^{p-1}}{t^{p-1}}\int_{-\infty}^\infty e^{-t\lambda_u(x)}\bigl(\lambda_u(x)-\Delta_u\bigr)^{-p}\lambda_u'(x)\ dx
\end{align*}
Using \eqref{E:A15} we thus get
\begin{multline}\label{E:A12}
\bigl\|\exp(-t\Delta_u)P_u\bigr\|_{0,2p} 
\\\leq\frac{\sqrt2\alpha_p(p-1)!}{2\pi}\frac{u^2}{t^{p-1}}e^{-\epsilon tu}\int_{-\infty}^\infty e^{-\epsilon tu|x|}\|(\lambda_u(x)-\Delta_u)^{-p}\|_{0,2p}\ dx
\end{multline}
From \eqref{E:A13} we obtain
\begin{equation}\label{E:A11}
\int_{-1}^1e^{-\epsilon tu|x|}\|(\lambda_u(x)-\Delta_u)^{-p}\|_{0,2p}\ dx
\leq 2C_p''(3\alpha_p+1)u^{p(p+1)}.
\end{equation} 
From \eqref{E:A10} we obtain
\begin{equation}\label{E:A9}
\int_{|x|\geq1}e^{-\epsilon tu|x|}\|(\lambda_u(x)-\Delta_u)^{-p}\|_{0,2p}\ dx
\leq\frac{2C_p'}{\epsilon}\frac1{tu}.
\end{equation}
Combining \eqref{E:A12}, \eqref{E:A11} and \eqref{E:A9} we find a constant $\bar C_p\geq0$ so that for all $u\geq u_0$ and $t>0$ we have
\begin{equation}\label{E:A8} 
\|\exp(-t\Delta_u)P_u\|_{0,2p}\leq\bar C_pe^{-\epsilon tu}\Bigl(\frac{u^{p(p+1)+2}}{t^{p-1}}+\frac u{t^p}\Bigr).
\end{equation}
Choose $\lambda_0\notin\Spec(\Delta_0)$, i.e.\ $\Delta_0-\lambda_0$ is invertible. Then, see \cite[Proposition~A.3.7]{S01},
$$
\|\exp(-t\Delta_u)P_u\|_{\tr}
\leq\|(\Delta_0-\lambda_0)^{-p}\|_{\tr}\|(\Delta_0-\lambda_0)^p\|_{2p,0}\|\exp(-t\Delta_u)P_u\|_{0,2p}.
$$
To complete the proof of Theorem~\ref{T:sesti} combine this with \eqref{E:A8} and note that $(\Delta_0-\lambda_0)^{-p}$ is trace class
since we assumed $p>n/2$.
\end{proof}

\section{Asymptotic of the large torsion and the proof of Theorem~\ref{T:main}}\label{S:assla}

Given the special role of the variable $u$ in this section we will replace the notation 
$\Delta_{q,u}$, $\tau_{\la,u}$ etc.\ by $\Delta_q(u)$, $\tau_{\la}(u)$ etc.

We will consider functions $t(u)$, $u\in(0,\infty)$ of the form
$$
\sum_{j=0}^nA_ju^j+\sum_{j=0}^nB_ju^j\log u+o(1),\qquad\text{as $u\to\infty$,}
$$
and refer to $A_0$ as the free term of $t(u)$ and denote it by $\FT(t(u))$. 
Note that $\log_\pi\tau_{\la}(u)$ is by Corollary~\ref{C:WHSla}
such a function. The key result of this section is Theorem~\ref{T:rella} below whose proof, although the same as of Theorem~B
in \cite{BFK96}, is supported by estimates derived in Section~\ref{S:whs}. This theorem calculates the free term
$$
\FT\bigl(\log_\pi\tau_{\la}(u)-\log_\pi\tilde\tau_{\la}(u)\bigr)
$$
provided  $M$ and $\tilde M$ have the same dimension, $E$ and $\tilde E$  the same rank, $f$ and $\tilde f$  the
same number of critical points in each index.

Here $\tau_{\la}(u)$ and $ \tilde\tau_{\la}(u)$ denote the large torsions associated with 
 two systems $(M,E,g,b,f)$ and $(\tilde M,\tilde E,\tilde g,\tilde b,\tilde f)$ as in Section~\ref{S:whs} 
which satisfy \eqref{E:M1}, \eqref{E:M2} and \eqref{E:M3}.

\subsection*{Asymptotic expansion of $\log\det$ for elliptic with parameter}

Suppose $E\to  M$ is a rank $k$ complex vector bundle over $(M,g)$ a smooth Riemannian manifold of dimension $n$,
$\mathbb D$ a second  order elliptic operator of Laplace--Beltrami type (cf.\ \cite{BFK96}) 
(i.e.\ the principal symbol $\sigma(\mathbb D)(\xi)=-\|\xi\|^2\id$),
$L\co E\to E$ a bundle map and $F\co M\to \mathbb R$ a smooth function. 
The operator $\mathbb D$ has  
always $\pi$ as a principal angle.\footnote{In fact any angle $\theta\neq0$ is a principal angle, i.e.\
for any $x$ and $\xi\neq0$ the spectrum of the finite dimensional linear map 
$\sigma(\mathbb D)(\xi,x):E_x\to E_x$ is disjoint from the ray of angle $\theta$.}
We denote also by $L$ resp.\ $F$ the zero order (differential)
operators defined by the bundle map $L$, resp. the multiplication by  $F$.

For any $u\in [0,\infty)$, let $\mathbb D (u):=\mathbb D^{L,F}(u)= \mathbb D  + u L + u^2 F$.  
If $F$ is strictly positive then the family $\mathbb D(u)$ is elliptic with parameter in the sense of \cite{S01}
or \cite{BFK92}.

We apply the considerations below to
$\mathbb D=\Delta+\epsilon$ where $\Delta:=\Delta_{E,g,b,q}$ is the $q$--Laplacian associated with 
the  flat  connection $\nabla^E$, the non-degenerate symmetric bilinear form $b$, and the Riemannian metric $g$  as defined in  
Section~\ref{S:intro}, and $\epsilon$  a positive real number.
The smooth function  $F$ will be $|df|^2$ where $f\co M\to\R$ is a smooth function on  
$M$ and the  endomorphism $L$ given by Lemma~\ref{L:F1}.  
In this case the family $\mathbb D (u)$  is elliptic with parameter away from the set of critical points of $f$.

With respect to a  coordinate  chart  $U\subseteq M$ with coordinates $x_1, \dotsc, x_n$ and a 
trivialization,  $E|_U= U\times \mathbb C^k$ the symbol  of the operator $\mathbb D$ is given by 
$$
\sigma(\mathbb D) (x,\xi)=-\|\xi\|^2  + \sum_{i=1}^n \alpha_i(x)\xi^i  + \beta(x),
$$  
where $\|\xi\|^2= \sum g^{ij}\xi^i\xi^j,$  $g^{ij}$ the Riemannian metric, $\alpha_i(x), \beta_i(x)$ are 
smooth  $\eend(\C^n)$--valued
functions.  The operator $\mathbb D(u)|_U$ is given by  
$$
\mathbb D(u)|_U=a_2(x,D,u) +a_1(x, D, u) +a_0(x,u),
$$
$x=(x_1,\dotsc,x_n)$, $D= (D_{x_1},\dotsc, D_{x_n})$,
where $a_{2-j}(x,\xi, u)$ , $\xi \in \mathbb R^n,$ are  $\eend(\mathbb C^k)$--valued smooth  functions with 
the homogeneity property 
$$
a_{2-j}(x,\tau\xi, \tau u)= \tau^{2-j} a_{2-j}(x,\xi, u),
$$
$\tau\in \mathbb R.$ Precisely   $a_2, a_1, a_0$ are given by 
\begin{equation}
\begin{aligned}a_2(x,\xi, u)= &-\|\xi\| ^2 + u^2 F\\
a_1(x,\xi,u)=& u L(x) +\sum_i^n\alpha_i(x)\xi^i\\
a_0(x,\xi, u)= &\epsilon +\beta(x).
\end{aligned}
\end{equation}

Suppose that $F$ is strictly positive, hence $\mathbb D (u)$ is elliptic with parameter.
As in \cite{BFK92} we define inductively the functions $r_{-2-j}(x,\xi,u,\mu)$,
$\mu\in\C$, with values in $\eend(\mathbb C^k)$ by:
\begin{align*}
r_{-2}(x,\xi, u, \mu)&:=(\mu-a_2(x,\xi,u))^{-1}\\
r_{-2-j}(x,\xi,u, \mu)&:=r_{-2}(x,\xi, u, \mu) \\
&\qquad\cdot\Biggl(\sum_{k=0}^{j-1}\sum_{|\alpha| +l +k= j}
\tfrac1{\alpha!}\partial^\alpha_{\xi} a_{2-l}(x,\xi,u,\mu) D^\alpha_x r_{-2-k}(x,\xi,u,\mu)\Biggr)
\end{align*}
with $\alpha$ a multi index $\alpha = (i_1,\dotsc,i_n)$.
Clearly $r_{-2-j}$ is homogeneous of degree $(-2-j)$ in $(\xi, u, \mu^{1/2})$.

We also define the smooth complex valued function 
\begin{equation}\label{E:E2}
a_{L,F}(x):= -\frac{1}{(2\pi)^n} \int _{\mathbb R^n} d\xi \int_0^\infty d\mu\tr\bigl(r_{-2-n}(x,\xi,1,-\mu)\bigr)
\end{equation} 
and as in \cite{BFK96}, page 352, equation (3.11) 
a simple calculation  shows   
\begin{equation}\label {E:E0}
a_{L,F}(x)+ a_{-L, F}(x)=0.
\end{equation}

Suppose  $M$ is closed and $(E, \mathbb D, F, L)$ as above. 
We can consider the complex valued function $\log\det_\pi\mathbb D(u)$.
The following result was established in the appendix of \cite{BFK92}, see also \cite{L97}.

\begin{theorem}[\cite{BFK92}]\label{T:101}
The  functions $a_{L,F}(x)$ define a density  on $M^n$.
If $\mathbb D(u)$ is invertible for $u$ large enough then 
$\log\det_\pi\mathbb D(u)$ has an asymptotic expansion of the form 
\begin{equation}\label{E:E1}
\sum _{j=0}^n A_j u^j + \sum_{j=0}^n B_j u^j \log u+ \sum_{i=1}^\infty C_i u^{-i},\qquad\text {as $u\to\infty$},
\end{equation}
with 
$A_0=\int_M a_{L,H}$.\footnote{Actually all terms $A_j$ and $B_j$ are integrals on $M$  of densities explicitly computable 
in each chart in terms of the symbol of $\mathbb D(u)$.} 
In particular $\log\det_\pi\mathbb D(u)$ is of the form 
\begin{equation}\label{E:E1'}
\sum _{j=0}^n A_j u^j + \sum_{j=0}^n B_j u^j \log u+o(1) \qquad\text {as $u\to\infty$}
\end{equation}
\end{theorem}
(The result proved in \cite{BFK92} is formulated under more general hypotheses and with stronger conclusions).
The result can be extended to compact manifolds with boundaries and Dirichlet boundary condition 
and  leads  to the following relative  version of Theorem~\ref{T:101} stated under more restrictive 
hypotheses  (satisfied in our situation).

\begin{theorem}[\cite{BFK96}]\label{T:102}
Suppose $(E_i,\mathbb D_i, L_i, F_i)$, $i=1,2$, are two systems as above. 
Suppose that there exist compact sets $K_i\subseteq M_i$ and open neighborhoods  $U_i$ of of $K_i$ so that:

1) $F_i |_{M\setminus K_i}$ are strictly positive, and

2) there exists the diffeomorphisms 
$ \varphi\co U_1 \to U_2$,  $\tilde\varphi\co E_1|_{U_1}\to E_2|_{U_2}$ 
bundle isomorphism above $ \varphi $ which intertwines 
$$
(E_1 |_{U_1}, \mathbb D_1,  L_1|_{U_1}, F_1|_{U_1})
\quad\text{with}\quad
(E_2 |_{U_2}, \mathbb D_2, L_2 |_{U_2}, F_2|_{U_2}).
$$
Let $V_i$ compact domains with smooth boundaries,  
$K_i \subseteq  V_i\setminus \partial V_i \subseteq V_i \subseteq U_i$ with $\varphi(V_1)= V_2$.
If\/ $\mathbb D_i(u)$ are invertible for $u$ large enough then  
$$
\log{\det}_\pi\mathbb D_1(u) - \log{\det}_\pi\mathbb D_2(u)
$$ 
has an asymptotic expansion of the form \eqref{E:E1} with 
$$
A_0= \int_{M_1\setminus V_1}  a_{L_1, F_1} - \int_{M_2\setminus V_2}  a_{L_2, F_2}
$$
\end{theorem}

\subsection*{A relative result for the asymptotic expansion of the large torsion}

Let $(M,E,g,b,f)$  be consisting of a closed smooth manifold $M$, complex flat bundle $E,$ Riemannian metric $g$, non-degenerate symmetric bilinear form $b$ and  Morse function $f$. We will refer to such  collection $(M,E,g,b,f)$ as a ``system'' provided
\eqref{E:M1}, \eqref{E:M2}, \eqref{E:M3} hold, and to a collection of neighborhoods $B_{q;j}$ of the critical points  
$x_{q;j}\in\mathcal X_q^f$ as $\rho$-admissible neighborhoods if they have disjoint closures and are the source of Morse charts 
$\varphi_{x_{q;j}}$ of radius $\rho$ so that on them \eqref{E:M1}, \eqref{E:M2} and \eqref{E:M3} hold.
Clearly such neighborhoods exists for $\rho$ small enough.

We define $B_{q;j}':=\varphi_{x_{q;j}}(D_{\rho/2})$, where $D_r$ denotes the disc of radius $r$ in $\R^n$ centered at $0$.

Given a collection of $\rho$--admissible neighborhoods $B_{q;j}$ introduce the manifolds
$$
M_I:=M\setminus\cup_{q,j}B'_{q;j};\quad M_{II}:=\cup_{q,j}\overline{B_{q;j}'},
$$
where $B'_{q;j}$ is defined as in the above definition. Both manifolds
$M_I$ and $M_{II}$ have the same boundary, given by a disjoint union of
spheres of dimension $n-1$.

Fix $\epsilon>0$ and consider the  operator $\Delta_q(u) +\epsilon$.
If $u$ is large enough, in view of Proposition~\ref{P:specgap} $\Delta_{q}(u)+ \epsilon$ 
has $\pi$ as an Agmon angle and is invertible. Its
symbol with respect to arbitrary coordinates $(\varphi,\psi)$ of
$(M,E\to M)$ is of the form
$$
a_2(x,\xi)+u^2\|\nabla f\|^2+a_1(x,\xi)+uL(x)+\epsilon
$$
where $a_i:B_{3\alpha}\times\mathbb R^d\to\eend(\Lambda^q(\mathbb R^d)\otimes\C^N)$,
$i=1,2$, are homogeneous of degree $i$ in $\xi$, where
$\|\nabla f\|^2:B_{3\alpha}\to\mathbb  R$ is given by
$$
\|\nabla f\|^2=\sum_{1\leq i,j\leq d}g^{ij}\frac{\partial f}{\partial x_i}
\frac{\partial f}{\partial x_j}
$$
and $L\co B_{2\alpha}\to\eend(\Lambda^q(\mathbb R^d)\otimes \mathbb C^N)$.

Therefore, away from the critical points of $f$, this operator is elliptic with parameter 
and  away from the critical points we can consider the densities 
$a_{L,H}$ associated with $(\Delta_q(u)+\epsilon)$; they will be denoted here by 
$a^q(f,\epsilon,x):=a^q(b,g,f,\epsilon,x).$ 
As in   \cite{BFK96} we can  establish the  following intermediary results:

\begin{proposition}\label{P:101} 
Assume that the systems $(M,E,b,g,f)$ and
$(\tilde M,\tilde E,\tilde b,\tilde g,  \tilde f)$ satisfy $\dim M= \dim \tilde M=n,$
$\sharp\mathcal X_q(f)=\sharp\mathcal X_q(\tilde f)$, $0\leq q\leq n$, and
$\rank E=\rank\tilde E$. Choose $\rho>0$ so that we have $\rho$--admissible neighborhoods of critical points for both systems.
Denote by $\Delta_q(u)$ resp.\ $\tilde \Delta_q(u)$ the Witten Laplacians associated with the two systems. 
Then, for any $\epsilon>0$, 
$$
\log{\det}_\pi(\Delta_q(u)+\epsilon)-\log{\det}_\pi (\tilde \Delta_q(u)+\epsilon)
$$ 
has an asymptotic expansion of the  form \eqref{E:E1} for $u\to\infty$ whose free term equals 
$$
\overline a^q(f,\tilde f,\epsilon)=\int_{M_{I}}a^q(f,\epsilon,x)-
\int_{\tilde M_{I}}a^q(\tilde f,\epsilon,\tilde x).
$$
\end{proposition}

For any $\epsilon>0$ introduce  
\begin{equation}
A(f,u, \epsilon):=\frac{1}{2}\sum_{q=0}^d(-1)^{q+1}q\log\det(\Delta_q(u)+\epsilon),
\end{equation}
which can be written as  $A(f,u,\epsilon)=A_\sm(f,u,\epsilon)+A_\la(f,u,\epsilon)$ with 
$$
A_{\sm/\la }(f,u,\epsilon):=\frac{1}{2}\sum_{q=0}^d(-1)^{q+1}q\log{\det}_\pi\bigl(\Delta_{\sm/\la,q}(u)+\epsilon\bigr)
$$
and $\Delta_{\sm/\la,q}(u):=\Delta_{q}(u) |_{\Omega^q_{\sm/ \la} (M; E_u)}$.
Clearly $A_\la(f,u,\epsilon=0)=\log\tau_{f,\la}(u)$.

\begin{proposition}\label{P:102}
With the assumptions in Proposition~\ref{P:101}, for any $\epsilon>0$ the quantities
$A(f,u,\epsilon)-A(\tilde f,u,\epsilon)$ and
$A_\la(f,u,\epsilon)-A_\la(\tilde f,u,\epsilon)$ have  
asymptotic expansions of the form \eqref{E:E1} for $u\to\infty$ which are identical. In particular  
$$
\FT\bigl(A(f,u,\epsilon)-A(\tilde f,u,\epsilon)\bigr)= 
\FT\bigl(A_\la(f,u,\epsilon)-A_\la(\tilde f,u,\epsilon)\bigr).
$$
\end{proposition}

\begin{proposition}\label {P:103}
With the assumptions in Proposition~\ref{P:102}:

(i) The limit
\begin{equation}
\lim_{\epsilon\to0}\FT\bigl(A_\la(f,u,\epsilon)-A_\la(\tilde f,u,\epsilon)\bigr)
\end{equation}
exists and is equal to $\FT\bigl(\log\tau_{f,\la}(u)-\log\tau_{\tilde f,\la}(u)\bigr)$.

(ii) This limit is given by 
$$
a(f,\tilde f)=\int_{M_{I}}\sum_q(-1)^qq a^q(f,\epsilon=0,x)-
\int_{\tilde M_{I}}\sum_q(-1)^qqa^q(\tilde f,\epsilon=0,\tilde x),
$$ 
( $\rho$  the same for both manifolds and small enough)
with $a^q(f,\epsilon,x)$ and $a^q(\tilde f,\epsilon,\tilde x)$ the densities considered above.
\end{proposition}

Combining the above propositions and (\ref{E:E0}) we have

\begin{theorem}\label{T:rella}
Assume that  the systems $(M^n,E,b,g,f)$ and $(\tilde M^n,\tilde E,\tilde b,\tilde g,\tilde f)$ satisfy
$\sharp\mathcal X_q(f)=\sharp\mathcal X_q(\tilde f)$, $0\leq q\leq n$, and
$\rank E=\rank\tilde E$.  Then 
$$
\FT\bigl(\log\tau_{\la}(u)-\log\tilde\tau_{\la}(u)\bigr)
$$
has an asymptotic expansion of the form \eqref{E:E1} whose free term is
$$
a(f,\tilde f)=\int_{M_{I}}\sum _q(-1)^qq a^q_0(f,\epsilon=0,x)-
\int_{\tilde M_{I}}\sum_q(-1)^qqa_0(\tilde f,\epsilon=0,\tilde x)
$$
with $a^q(f,\epsilon,x)$ and $a^q(\tilde f,\epsilon,\tilde x)$  the densities considered 
in Proposition~\ref{P:101} above.

(iii) When  $\dim M=n$ is odd we have
$$
a(f,\tilde f)+ a(n-f,n-\tilde f)=0.
$$
\end{theorem}

Proposition~\ref{P:101} is similar to Proposition 3.1 in \cite{BFK96} (but for non-selfadjoint Witten Laplacians)
so the proof is the same. It is actually  a straightforward consequence of the Theorem~\ref{T:102} above.

The proof of Proposition~\ref{P:102} goes as follows. 
As the eigenvalues of the operator $\Delta_{\sm,q}(u)$ tend exponentially fast to $0$ as
$u\to\infty$ in view of Proposition~\ref{P:specgap}  
$$
\log_\pi \det(\Delta_{q,\sm}(u)+\epsilon)=m_q\log\epsilon+O\left(\frac{1}{\epsilon}\alpha e^{-\beta u}\right)
$$
for some positive constants $\alpha, \beta$, and therefore,
$A_\sm(f,u,\epsilon)-A_\sm(\tilde f,u,\epsilon)$ is exponentially small as
$u\to\infty$.  Therefore for any $\epsilon>0$,
$$
A(f,u,\epsilon)-A(\tilde f,u,\epsilon)
$$ 
and
$$
A_\la(f,u,\epsilon)-A_\la(\tilde f,u,\epsilon)
$$ 
have asymptotic expansions of the form \eqref{E:E1} for $u\to\infty$ which
are identical. q.e.d.

Proposition~\ref{P:103} is more elaborated and the remaining of the subsection elaborate on this proof.

\begin{proof} [Proof of Proposition~\ref{P:103}]
To check (i) we verify that the function $H(u,\epsilon)$, defined for $\epsilon>0$ and $u$
sufficiently large by
$$
H(u,\epsilon):=A_\la(f,u,\epsilon)-A_\la(\tilde f,u,\epsilon)+\log \tau_{\la}(u)
-\log \tilde \tau_{\la}(u)
$$
is of the form
\begin{equation}\label{E:206}
H(u,\epsilon)=\sum_{k=1}^T\epsilon^kf_k(u)+g(u,\epsilon),
\end{equation}
where $g(u,\epsilon)=O( u^{-1+\delta})$ uniformly in $\epsilon$. The statement of
Proposition~\ref{P:103} can be deduced from  this formula  as follows: Recall
that for $\epsilon>0$, $H(u,\epsilon)$ has an asymptotic expansion for
$u\to\infty$ because $A_\la(f,u,\epsilon)-A_\la(\tilde f,u,\epsilon)$ and $\log \tau_{\la}(u)$
and $\log \tilde \tau_{\la}(u)$ have, the first by Theorem~\ref{T:rella} the last two by Corollary~\ref{C:WHSla}.
As $g(u,\epsilon )=O(u^ {-1+\delta})$ uniformly in $\epsilon$ we
conclude that for any $\epsilon>0$, $\sum_{k=1}^T\epsilon^kf_k(u)$ has an asymptotic
expansion for $u\to\infty$. By taking $T$ different values
$0<\epsilon_1<\dots<\epsilon_T$ for $\epsilon$ and using that the Vandermonde
determinant is nonzero
$$
\det \left( \begin{smallmatrix}
\epsilon_1 & \dots & \epsilon_1^T\\
\vdots && \vdots\\
\epsilon_d & \dots & \epsilon_d^T
\end{smallmatrix}\right)
\neq 0
$$
we conclude that for any $1\leq k\leq T$, $f_k(u)$ has an asymptotic
expansion for $u\to\infty$ and that for any $\epsilon>0$
$$
\FT(H(u,\epsilon))=\sum_{k=1}^d\epsilon^k\FT(f_k(u)).
$$
Hence $\lim_{\epsilon\to0}\FT(H(u,\epsilon))$ exists and
$\lim_{\epsilon\to0}\FT(H(u,\epsilon))=0$.
It remains to prove \eqref{E:206}.

Recall that if
$$
\theta_{q,\la}(u,\mu):=\tr(e^{-\mu\Delta_q(u)} P_u)
$$
and 
\begin{equation}\label{E:201}
\zeta_{q,\la}(u,\epsilon,s):=\frac1{\Gamma(s)}\int_0^\infty \mu^{s -1}\theta_{q,\la}(u,\mu)  e^{-\epsilon \mu} d\mu
\end{equation}
then  
\begin{equation}\label {E:202}
\log{\det}_\pi(\Delta_{q, \la}(u)+\epsilon)= \frac {d}{ds}\Big|_{s=0}\zeta_{q,\la}(u,\epsilon, s).
\end{equation}
We have the estimate

\begin{lemma}\label{L:201}
There exists a  constant $C_1,C_2,\beta,>0$, $0<\delta <1$, and integer $T\geq n$ so that

(i) for $u$ large enough  and $0\leq \mu \leq u^{-1+\delta}$
\begin{equation}\label{E:203}
\theta_q(u,\mu)\leq C_1\mu^{-T}
\end{equation}

(ii)  for $u$ large enough, and $\mu \geq u^{-1+\delta}$
\begin{equation}\label{E:204}
\theta_q(u,\mu)\leq C_2e^{-\beta u\mu}.
\end{equation}
\end{lemma}

\begin{proof}
Note that Theorem~\ref{T:sesti} implies  that  there exists an integer $N \geq n$
and the constants $C', \beta >0$ so that for $u$ large enough 
$$ 
| \theta (u ,\mu)| \leq C' e^{-\gamma u \mu}  u^N/\mu^n
$$ 
with $\gamma >0$. Choose $0<\delta<1$, (for example $\delta= 1/2)$.

Suppose $\mu \leq u^{-1 +\delta}$.  This implies $\mu u^{1-\delta} \leq 1$; 
and supposing  $u$ large enough  we have  $\mu \leq 1/2$.

As  $\mu\leq u^{\delta-1}$, hence $ u^N\leq \mu^{-N/(1-\delta)}$,  we get 
$$
|\theta(u,\mu)|\leq C'e^{-\gamma \mu u} /  u^{n+ N/(1-\delta)} \leq  C' e / u^T,
$$
where $T> n+ N/(1-\delta)$. This establishes (i).

Suppose $\mu\geq u^{-1 +\delta}$. This implies $\mu^{-n} \leq u^{n(1-\delta)}$ and therefore 
$$
e^{-\gamma/2 \  \mu  u}u^N/\mu^n \leq e^{-\gamma/2\  u^\delta}  u^{N+n(1-\delta)}.
$$
Clearly for $u$ large enough  we have  
$$
e^{-\gamma/2 \  u^\delta}  u^{N+n-\delta} \leq 1
$$ 
since $\delta > 0$. Take $C_2= C'$, $\beta= \gamma/2$  and (ii) is verified.
\end{proof}

To establish  \eqref{E:206} we decompose the function $\zeta_{q,\la}(u,\epsilon,s)$ into two parts
\begin{equation}\label{E:207}
\zeta_{q,\la}^I(u,\epsilon,s)=\frac{1}{\Gamma(s)}\int_{u^{-1+\delta}}^\infty  \mu^{s-1}\theta_q(u,\mu)
e^{-\epsilon\mu}d\mu
\end{equation}
and
\begin{equation} \label{E:208}
\zeta_{q,\la}^{II}(u,\epsilon,s)=\frac{1}{\Gamma(s)}\int_0^{u^{-1+\delta}}\mu^{s-1}\theta_q(u,\mu)
e^{-\epsilon\mu}d\mu.
\end{equation}

First let us consider
$$
\zeta_{q,\la}^I(u,\epsilon,s)-\zeta_{q,\la}^I(u,\epsilon=0,s)=
\frac{1}{\Gamma(s)}\int_{u^{-1+\delta}}^\infty\mu^s\theta_q(u,\mu)
\frac{e^{-\epsilon\mu}-1}{\mu}d\mu
$$
Note that
$$
\zeta_{q,\la}^I(u,\epsilon,s)-\zeta_{q,\la}^I(u,\epsilon=0,s)
$$
is by Lemma~\ref{L:201}(ii) an entire function of $s$ and so is   $1/\Gamma(s)$.

Clearly, $\tfrac1{\Gamma(s)}\big|_{s=0}=0$,
$\tfrac{d}{ds}\big|_{s=0}\tfrac1{\Gamma(s)}=1$
and $1-e^{-\epsilon\mu}\leq\epsilon\mu$.
Therefore by Lemma~\ref{L:201} we have
\begin{multline*}
\Biggl|\frac{d}{ds}\Big|_{s=0}
\bigl(\zeta_{q,\la}^I(u,\epsilon,s)-\zeta_{q,\la}^I(u,\epsilon=0,s)\bigr)\Biggr|
\\
=\Bigl|\int_{ u^{-1+\delta}}^\infty\theta_q(u,\mu)\frac{e^{-\epsilon \mu}-1}{\mu}d\mu\Bigr|
\leq\epsilon C_2 \int_{ u^{-1 +\delta}}^\infty e^{-\beta u\mu}d\mu
=\frac{\epsilon C_2}{\beta u}e^{-\beta u^\delta}.
\end{multline*}

Concerning the term
$$
\frac{d}{ds}\Big|_{s=0}\big(\zeta_{q,\la}^{II}(u,\epsilon,s)-\zeta_{q,\la}^{II}(u,\epsilon=0,s)\bigr),
$$
expand $(e^{-\epsilon\mu}-1)/\mu$
$$
(e^{-\epsilon\mu}-1)/\mu=\sum_{k=1}^T\frac{(-1)^k}{k!}\epsilon^k\mu^{k-1}
+\epsilon^{T+1}\mu^de(\epsilon,\mu)
$$
where the error term is given by
$$
e(\epsilon,\mu)=\left(\sum_{k=T+1}^\infty\frac{(-1)^k}{k!}\epsilon^k\mu^{k-1}\right)
/\epsilon^{T+1}\mu^T.
$$
Note that  according to Lemma~\ref{L:201}(i) we have $\mu^T \theta_q(u,\mu) \leq C_1$.

Therefore
$$
\int_0^{ u^{-1+\delta}}\mu^s\theta_q(u,\mu)\epsilon^{d+1}\mu^Te(\epsilon,\mu)d\mu
$$
is a meromorphic function of $s$, with $s=0$ a regular point and, for  sufficiently large $u$ we have 
$$
\Biggl|\frac{d}{ds}\Big|_{s=0}\left(\frac{1}{\Gamma(s)}
\int_0^{ u^{-1+\delta}}\mu^s\theta_q(u,\mu)
\epsilon^{T+1}\mu^Te(\epsilon,\mu)d\mu\right)\Biggr|
\leq\epsilon ^{T+1}C_1 u^{-1+\delta}
$$
with $C_1$ is independent of $u$ and $\epsilon $, $0\leq\epsilon\leq1$.

Finally, recall that $\theta_q(u,\mu)$ admits an expansion for $\mu\to 0+$ of the form
$$
\theta_q(u,\mu)=\sum_{j=0}^TC_j(u)\mu^{(j-d)/2}+\theta_q'(u,\mu)
$$
where $\theta_q'(u,\mu)$ is continuous in $\mu\geq0$. Therefore, for $1\leq k\leq T$,
$$
\frac{1}{\Gamma(s)}\int_0^{u^{-1+\delta}}\mu^s\theta_q(u,\mu)\frac{(-1)^k}{k!}\epsilon^k\mu^{k-1}d\mu
$$
is analytic with respect to $s$ at $s=0$ and
$$
\sum_{k=1}^T\frac{d}{ds}\Big|_{s=0}\Biggl(\frac{1}{\Gamma(s)}
\int_0^{u^{-1+\delta}}\mu^s\theta_q(u,\mu)\frac{(-1)^k}{k!}\epsilon^k\mu^{k-1}d\mu\Biggr)
$$
is of the form $\sum_{k=1}^T\epsilon ^kf_k(u)$. This establishes \eqref{E:206}.
Part (ii) follows from Proposition~\ref{P:101} and  part (iii) from (\ref{E:E0}).
\end{proof}

\subsection*{Proof of Theorem~\ref{T:main}}

In this section we want to check that $S_{E,[b]}^2=1$.
Consider a system $(M, E, g, b, f)$.  Clearly $(M,E,g,b,-f)$ is also a system  and denote by 
$\tau^+_{\la}(u)$ resp $\tau^-_{\la}(u)$ the large torsion for the first resp.\ of the second. Similarly we write 
$\tau(\Int^+_{\sm,u})$ resp.\ $\tau(\Int^-_{\sm,u})$ for the relative torsion of \eqref{E:intusm} in Section~\ref{S:intro} 
when applied to the first resp.\ second system.

Suppose now we have two systems  $(M,E,g,b,f)$ and $(\tilde M,\tilde E,\tilde g,\tilde b,\tilde f)$, 
and suppose that $f$ and $\tilde f$ have the same number of critical points in each index. Since $\dim M=n$ is odd we have
$(-X)^*\Psi_g=-X^*\Psi_g$, and therefore Corollary~\ref{C:WHSla} yields
\begin{equation}\label{E:K11}
S^2_{E,[b]}/S^2_{\tilde E,[\tilde b]}
=\frac{\tau^+_{\la}(u)\cdot\tau^-_{\la}(u)}{\tilde\tau^+_{\la}(u)\cdot\tilde\tau^-_{\la}(u)}\cdot e^{\beta u}\cdot\bigl(1+ O(e^{-\epsilon u})\bigr)
\end{equation}
with a real number $\beta$.
We know that  the  left  side of \eqref{E:K11} is constant and
$$
\log_\pi\frac{\tau^+_\la(u)\cdot\tau^-_\la(u)}{\tilde\tau^+_\la(u)\cdot\tilde\tau^-_\la(u)}
$$ 
has an asymptotic expansion whose free part is by Theorem~\ref{T:rella} equal to $0$.
This implies that 
$$
\mathcal S_{E,[b]}^2/S_{\tilde E,[\tilde b]}^2=1.
$$

We choose $\tilde M$ to be the sphere, $\tilde E$ the trivial flat bundle of the same rank as $E$ and $\tilde b$ the canonical 
symmetric bilinear form. We note that since $\dim M=n$ is odd, one can always provide a Morse function $\tilde f$ on $\tilde M:=S^{\dim M}$ with the same 
number of critical points as $f$ in each index. Since $\mathcal S^2_{\tilde E,[\tilde b]}=1$ we conclude that $\mathcal S_{E,[b]}^2=1$.
q.e.d.

\section {Appendix; Remarks on the proof of Conjecture~\ref{C:main}}

It is likely that in Theorem~\ref{T:main} one can replace $\pm 1$ by $1$. For this purpose  
notice that for each Morse function $f$ one can produce a  square root 
$S'_{E, [b],f}$ of $S_{E [b]}$  and by the same arguments as in the proof of Theorem~\ref{T:main} one can show that 
$$
S'_{E,[b],f} \cdot S'_{E,[b], -f}=1.
$$
Clearly, if we show that $S'_{E, [b], f}$ is independent of $f,$ then $S_{E,[b]}=1$.

To define a square root of $S_{E,[b]}$ we use the formulas \eqref{E:intusm1} and \eqref {E:DD}. We need 
two additional data: a Morse--Smale vector field $X$ and an orientation of the total space of  
the  mapping cone complex associated with the quasi isomorphism  $\Int_{\sm,u}\co\Omega_{\sm}(M;E_u)\to C(X;E_u)$. 
The  orientation  provides a square root of $\tau(\Int_{\sm,u})$. Both $\tau_{\la,u}$ and 
$\exp\bigl(-2\int_{M\setminus \mathcal X} \omega_{E_u,b}\wedge (-X)^\ast \Psi_g\bigr)$ 
have an unambiguous square root. 
We can choose $X=-\grad_g f$ in which case  $\Int_{\sm,u}$ is an isomorphism for $u$ large enough and a canonical  
orientation is implicit in the construction of the mapping cone. The product of these square roots give our desired square root.

It is possible to show that this construction extends to generalized Morse  functions.
Recall that  by a result of H.~Chaltin \cite{HKW95} any two Morse functions $f_1$ and $f_2$ can be 
joined by a homotopy $f_t$ with $f_t$  Morse function for all $t\in [1,2]$ but $t_1,t_2,\dotsc,t_k$, and
generalized Morse function for $t=t_1,t_2,\dotsc,t_k$. One can even arrange that each generalized Morse 
function has only one birth/death or death/birth critical point. Then the independence of $S'(E,[b],f)$ of 
the  Morse function  $f$ follows  for from the continuity in $t$ of $S'_{E,[b], f_t}$ for a homotopy of the 
type provided by Chatlin result.

Recall that a generalized Morse function  is  a  smooth function whose critical points are either non-degenerate or 
are birth-death/death-birth critical points, cf.\ \cite{HKW95}.  For each $q$,  let $m_q$ be the number of non-degenerate 
critical points of index $q$ and $m'_q$ the number of degenerate (birth/death) critical points  of index $q$.
It was established in \cite{HKW95}  that for 
$u$ large the spectrum of $\Delta_{E_u,g,b,q}$ decomposes in three disjoint parts  $\Spec_{\sm,u}$,  $\Spec_{\rm{mla}, u}$,  
$\Spec_{\rm{vla}, u}$, called small spectrum, moderately large  spectrum and and very large spectrum.  The small spectrum   
$\Spec_{\sm, u}$ consists of  $(\rank E)\cdot   m_q$ complex numbers converging exponentially fast to $0$,  
$\Spec_{\sm, u}$ consists of $2 (\rank E)\cdot  m'_q$ complex numbers whose both real part and absolute value are 
converging to $\infty$ but  not faster than $Cu^{2/3}$ and $\Spec_{\rm{vla}, u}$, the  rest of the spectrum,  
consists  of complex numbers  whose real part is  converging to $\infty$ faster than $C'  u,$ with $C, C'$ some constants. Actually this was established  for a flat vector bundle equipped with a Hermitian structure 
but  we expect  the statements  hold true  for a  non-degenerate symmetric bilinear form 
also.

To define the square root of $S_{E,[b]}$ for such generalized  Morse function one can use 
instead of $\Omega_{\sm}(M;E_u)$ and $\Omega_{\la}(M;E_u)$ either $\Omega_{\sm}(M;E_u)$ and 
$\Omega_{\rm{mla}}(M;E_u)\oplus \Omega_{\rm{vla}}(M;E_u)$ or $\Omega_{\sm}(M;E_u)\oplus 
\Omega_{\rm{mla}}(M;E_u)$ and $\Omega_{\rm{vla}}(M;E_u)$.
One can choose conveniently a Morse--Smale vector field, for example a vector field which away from the 
degenerate critical points is gradient like for the generalized Morse function; in this case 
canonical orientations exist for the mapping cone of the corresponding integration morphisms   in both cases. 
The square roots obtained by either choice are the same.

While the continuity at $t'$ when $f_{t'}$ is a Morse function is straightforward at $t'$ when $f_{t'}$ is 
generalized Morse function is more subtle. It is easier to check this continuity when there is only one  
birth/death or dearth/birth critical point and this can be done  separately from left and from right using 
the two possible definitions of the  square root.

An extension of the result $S_{E,[b]}=1$ to smooth  odd dimensional manifolds  with boundary 
combined with a product formula  for $S_{E,b},$ (like the product formula for for analytic or 
combinatorial  torsion,)  will imply the result for even dimensional manifolds as well. 
The details of the above remarks will be presented  in a forthcoming paper.

\end{document}